\documentclass{article}

\usepackage{PRIMEarxiv} % Template used

\usepackage[utf8]{inputenc} % allow utf-8 input
\usepackage[T1]{fontenc}    % use 8-bit T1 fonts
\usepackage[hypertexnames=false,colorlinks=true,linkcolor=blue,citecolor=blue]{hyperref} % blue hyperrefs
\usepackage{url}
\usepackage{booktabs}
\usepackage{amsfonts}
\usepackage{nicefrac}
\usepackage{microtype}
\usepackage{lipsum}
\usepackage{fancyhdr}       
\usepackage{graphicx}      

\usepackage{subfloat}
\usepackage{svg}
\usepackage{tikz}

\usepackage{lmodern}
\usepackage{iftex}
\usepackage[english]{babel}
\usepackage[utf8]{inputenc}
\usepackage{float}
\usepackage{caption}
\usepackage{subcaption}
\usepackage{listings}
\usepackage{gensymb}
\usepackage[labelfont=bf]{caption}
\usepackage{ragged2e}
\usepackage{booktabs}
\usepackage{makecell}
\usepackage{fancyhdr}
\usepackage{trimclip}
\usepackage{color, colortbl}
\usepackage[numbers, sort&compress]{natbib}
\usepackage{epsfig}
\usepackage{makecell}
\usepackage{amssymb}
\usepackage{amsmath}
\usepackage{xcolor}

\DeclareMathOperator*{\argmin}{arg\,min}

\usepackage{algorithm}
\usepackage{algpseudocode}

\lstdefinestyle{code}{
language={Matlab},
basicstyle=\ttfamily\footnotesize,
keywordstyle=\color{blue}
}
\usetikzlibrary{shapes.geometric, arrows}

\tikzstyle{normal} = [rectangle, rounded corners, 
minimum width=2cm, 
minimum height=1cm,
text centered, 
draw=black, 
fill=blue!30]

\tikzstyle{small} = [rectangle, rounded corners, 
minimum width=2cm, 
minimum height=0.5cm,
text centered, 
draw=black, 
fill=blue!20]

\tikzstyle{init} = [rectangle, minimum width=3cm, minimum height=1cm, text centered, draw=black, fill=orange!30]

\tikzstyle{decision} = [rectangle, minimum width=1.5cm, minimum height=1cm, text centered, draw=black, fill=green!30]

\tikzstyle{arrow} = [thick,->,>=stealth]
\tikzstyle{line} = [thick, -]
\usetikzlibrary{calc}

\title{Data-Driven Stabilisation of Unstable Periodic Orbits of the Three-Body Problem}

\begin{document}
\maketitle
\vspace{-2.5cm}
\begin{center}
 \textbf{Owen M. Brook$^{1}$, Jason J. Bramburger$^{2}$, Davide Amato$^{1}$, Urban Fasel$^{1}$} \\
  { \fontsize{9}{10}\selectfont  \textit{$^1$ Department of Aeronautics, Imperial College London, London, UK \\ $^2$ Department of Mathematics and Statistics, Concordia University, Montr\'{e}al, QC, Canada} } \\
  { \fontsize{9}{10}\textit{Corresponding author: Owen Brook (owen.brook20@imperial.ac.uk)}} \\
\end{center}

\noindent

\vspace{0.5cm}
\begin{abstract}
Many different models of the physical world exhibit chaotic dynamics, from fluid flows and chemical reactions to celestial mechanics. The study of the three-body problem (3BP) and its many different families of unstable periodic orbits (UPOs) have provided fundamental insight into chaotic dynamics as far back as the 19th century. The 3BP, a conservative system, is inherently challenging to sample due to its volume-preservation property.
In this paper we present an interpretable data-driven approach for the state-dependent control of UPOs in the 3BP, through leveraging the inherent sensitivity of chaos and the local manifold structure. We overcome the sampling challenge by utilising prior knowledge of UPOs and a novel augmentation strategy.
This enables sample-efficient discovery of a verifiable and accurate local Poincar\'e map in as few as 55 data points. We suggest that the Poincar\'e map is best discovered at a surface of section where the norm of the monodromy matrix, i.e. the local sensitivity to small perturbations, is the smallest.
To stabilise the UPOs, we apply small velocity impulses once each period, determined by solving a convex system of linear matrix inequalities based on the linearised map. We constrain the norm of the decision variables used to solve this system, resulting in locally optimal velocity impulses directed along the local stable manifold. Critically, this behaviour is achieved in a computationally efficient manner.
We demonstrate this sample-efficient and low-energy method across several orbit families in the 3BP, with potential applications ranging from robotics and spacecraft control to fluid dynamics.

\end{abstract}

\section{Introduction}
\label{sec: Intro}
Chaotic systems are characterised by their sensitivity to small perturbations, where two trajectories starting arbitrarily close together will diverge exponentially in time. For a long time, this meant that these systems were largely avoided and efforts were generally directed to suppress or evade chaos \cite{Rega2010ControllingDynamics}. However, chaotic systems are unavoidable and appear all throughout nature, from fluid dynamics \cite{Fazendeiro2010UnstableTurbulence}, and the movement of comets \cite{Koon2000HeteroclinicMechanics}, through to the dynamics of the heart \cite{Wagner1998ChaosUpdate}, and chemical reactions \cite{Cramer2006ChaosReview}. 
During the 19th century, the pioneering work of Henry Poincar\'{e} studied such systems by considering the behaviour of different sets of initial conditions, rather than just an individual trajectory. Considering the flow of an $N$ dimensional system, Poincar\'e introduced the idea of a surface of section, a $N - 1$ dimensional hypersurface perpendicular to the flow, and studying a discrete map which describes successive crossings of this surface. Fixed points of this so called Poincar\'e map were shown to be either equilibria or periodic orbits of the original system, and through this, much of our knowledge of chaos has emerged \cite{Ott2002ChaosSystems}. A key understanding is that chaotic systems are made up of infinite numbers of unstable periodic orbits (UPOs), which have been shown to drive chaotic dynamics in turbulent flow \cite{Budanur2017RelativeFlow}, convection \cite{Singer1991ActiveConvection}, and the three-body problem which Poincar\'{e} originally studied.

\textbf{Chaos control.} Building on Poincar\'e's research nearly a century earlier, the work of Ott, Grebogi, and Yorke, harnessed the sensitivity of chaotic systems, using small perturbations to control them \cite{Ott1990ControllingChaos}. That is, if a small perturbation causes a large change in response, then we only need to apply a small perturbation to control such a system. This method, often abbreviated as the OGY method, consists of applying small parameter perturbations each time a trajectory crosses a surface of section within a neighbourhood of the desired UPO's states.
The OGY method and its various extensions have successfully demonstrated control of many real-world chaotic systems, such as coupled mechanical systems \cite{Boccaletti2000TheApplications}, a rabbit heart \cite{Weiss1994ChaosBiology.}, and atomic force microscopy \cite{Arjmand2008ChaosControl}. Different approaches have also been proposed with different goals in mind, such as eliminating or shifting bifurcations \cite{Rega2010ControllingDynamics}.
Through exploiting this sensitivity and using Poincar\'e maps, the OGY method results in significantly reduced control effort compared to traditional control methods, beneficial in many engineering applications. 

\textbf{Challenges.} In spite of this success, wider applications of OGY-style control methods remain a challenge as they require an analytical expression of a Poincar\'e map, which cannot be derived other than for the most simple systems. Recently, various data-driven techniques have been used to learn these maps, for example reinforcement learning was used to learn the mappings of the walking pattern of a bipedal robot and control it \cite{Morimoto2005Poincare-Map-BasedWalking}, and a novel neural network architecture, H\'{e}nonNet, used to learn maps of toroidal magnetic fields \cite{Burby2020FastFields}, as well as other deep learning methods \cite{Bramburger2021DeepMappings}. Recent work employed the dictionary-learning model identification method, \emph{sparse identification of nonlinear dynamics} (SINDy) \cite{Brunton2015DiscoveringSystems}, to discover a Poincar\'{e} map from data \cite{Bramburger2020PoincareTheory}. This enabled the identification of UPOs and their subsequent control \cite{Bramburger2020Data-DrivenOrbits} using the pole-placement method \cite{Romeiras1992ControllingSystems} and linear matrix inequalities \cite{Parrilo2012SemidefiniteOptimization} (see \cite{Bramburger2024IdentifyingDynamics} for a full overview). Using SINDy to discover Poincar\'e maps has potential benefits over deep learning methods, due to its interpretability, the low amount of data required and the robustness of various versions of the algorithm \cite{Fasel2021Ensemble-SINDy:Control, Kaheman2020SINDy-PI:Dynamics}. As such, SINDy is better suited to settings where collecting data comes at a significant cost either experimentally or computationally, and where noise is present.
The work of \cite{Bramburger2020Data-DrivenOrbits} draws inspiration from Ott, Grebogi, and Yorke's original work and almost exclusively focuses on their application to parameter-dependent dissipative chaotic systems such as the R\"ossler and a chaotic jerk system. 
In many real-world systems, only the system's states, not their parameters, can be altered and so further investigating the state-dependent control of UPOs is of interest.
Learning maps and controlling Hamiltonian and conservative systems also remains a significant challenge due to their volume-preserving nature. These systems are less forgiving, because unlike dissipative systems, random initial conditions left to evolve will not all fall onto a regular attractor structure. Examples of such systems include the ray equations, magnetic field lines in plasma, and characteristics of statistical mechanics~\cite{Ott2002ChaosSystems}.

\textbf{Control of the three-body problem.} The final example in \cite{Bramburger2020Data-DrivenOrbits} applies their data-driven OGY method to a problem of station keeping wherein a stroboscopic (periodic) mapping was learned and precise periodic impulses were applied to force a spacecraft to rest at an equilibrium point of the three-body problem (3BP). This conservative system describes the motion of three particles under the influence of attractive forces from each other. It admits five equilibria known as Lagrange, or libration, points and contains many different families of periodic orbits, both stable and unstable. Analogues of the 3BP appear not just in celestial mechanics but in many different fields, from quantum mechanics \cite{Aquilanti2004Three-bodySets} to protein folding \cite{Ejtehadi2004Three-bodyModels}, with new families of orbits still being found \cite{Dzhanoev2010ChaosProblem} along with methods for computing them \cite{Baresi2020AProblems}. In a somewhat similar approach to the OGY method, the chaotic properties of the 3BP have been leveraged to plan missions which take spacecraft along chaotic trajectories \cite{Marsden2006NewDesign}, using only a small number of velocity impulses. The existence of chaotic trajectories connecting UPOs of the 3BP, known as heteroclinic connections, were proven in \cite{Koon2000HeteroclinicMechanics}, allowing spacecraft to jump between different orbit families, such as in the NASA Genesis mission which had a fuel mass percentage of only five percent \cite{Ross2006TheNetwork}. 
% Although these trajectories take significantly longer to transit along, the savings from fuel expenditures can offer an advantageous trade-off. 
Keeping spacecraft along UPOs forms a critical part of such missions, and there are a plethora of methods to achieve this \cite{Shirobokov2017SurveyOrbits}. Floquet mode control utilises knowledge of underlying dynamics to find an impulse which ensures the component of the error vector after one period is zero \cite{Simo1987OnOrbits}. Another alternative method linearises the dynamics at various target points close to the reference UPO and computes impulses that minimise the deviation from the UPO at these points \cite{Howell1993StationkeepingTrajectories}.
The application of these methods to actual missions is complicated by the sensitivity of the dynamics to small perturbations, manifesting in short timescales for divergence from reference UPOs \cite{Folta2014EarthMoonOperations}. Consequently, classical station-keeping methods are hard to generalise to more sophisticated models of the 3BP including primaries on eccentric orbits and additional perturbations, leading to suboptimal solutions. Other methods often rely on computationally expensive direct optimisation approaches to find an optimal impulse \cite{Pavlak2012StrategySystem}, and are not easily generalised to other chaotic systems.

\textbf{Our contributions.} In this paper we present a novel, sample-efficient and interpretable data-driven approach that achieves state-dependent control of UPOs of the conservative 3BP. Our method utilises as few as 55 data points to discover an accurate and interpretable Poincar\'e map quickly by using prior knowledge of UPOs and a novel augmentation strategy, offering a key advantage over black-box deep learning methods. We find that a surface of section on which the Poincar\'e map is discovered is best placed where the sensitivity to perturbations is the least, which we quantify by taking the norm of the monodromy matrix. We exploit the information about the underlying dynamics contained in this map to perform low-energy locally-optimal velocity impulses directed along the local stable manifold, by constraining the decision variables of a series of linear matrix inequalities. We demonstrate our method by stabilising both planar Lyapunov and non-planar halo UPOs in the 3BP, one of the most widely studied conservative and Hamiltonian chaotic systems. Our intuitive method can readily be applied to more sophisticated models of multi-body systems and other chaotic volume-preserving and dissipative systems.
% displays many of the key features of more complex/realistic models of celestial mechanics.
\begin{figure}
    \centering
    \includegraphics[width=1.0\linewidth]{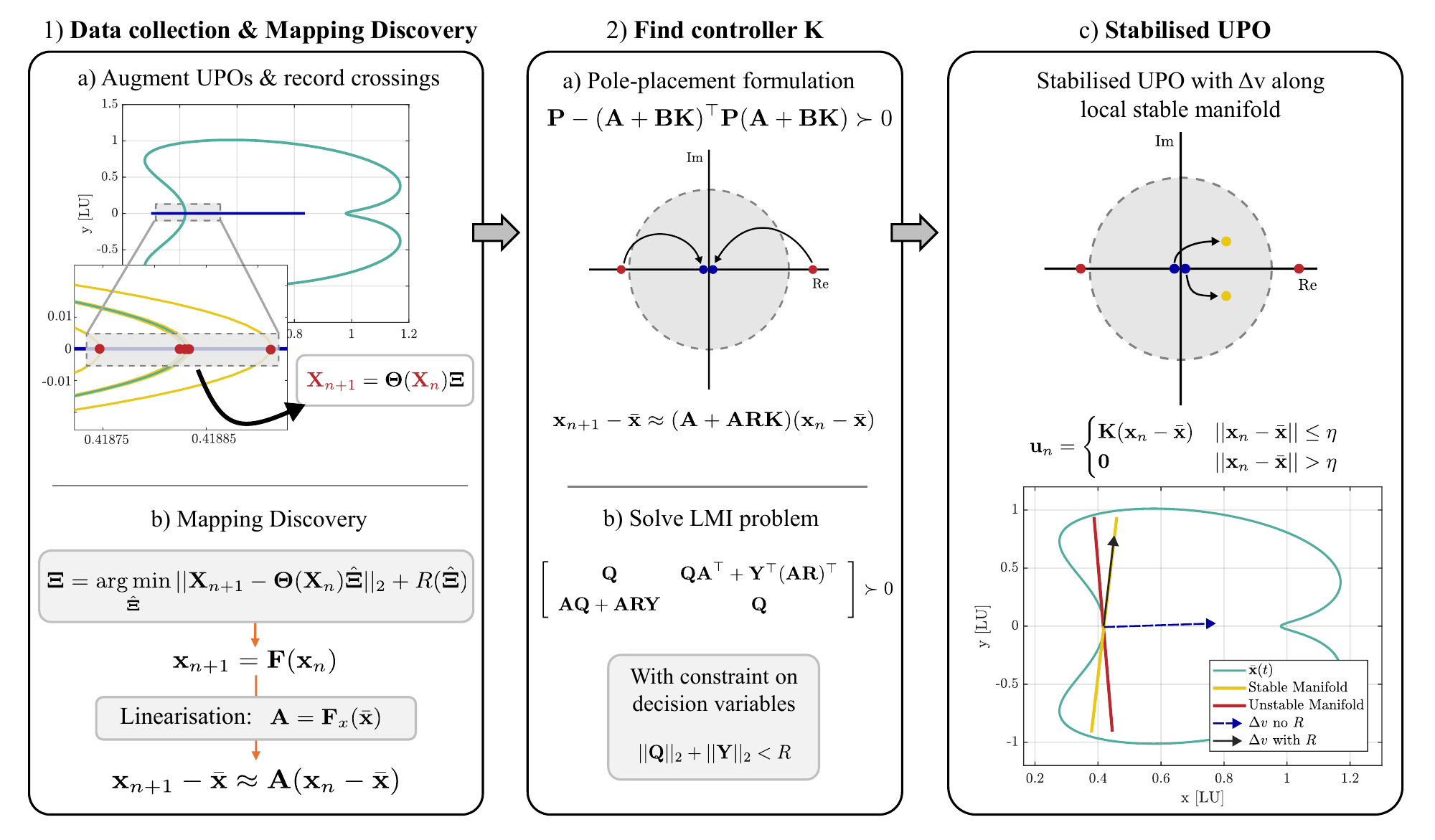}
    \caption{Our method for stabilising a Lyapunov UPO: 1) Initial conditions of known UPOs are selected and small velocity perturbations added. The states $(x,\Dot{x},\Dot{y})$ are recorded at successive crossings of the surface of section $\Sigma =  \{ (x,\Dot{x},\Dot{y})| y = 0,  x < 0.84\}$. A Poincar\'e map is discovered using SINDy and linearised at the desired UPO $\mathbf{\bar{x}}$. 2) A controller $\mathbf{K}$ is found through solving a pole-placement problem using linear matrix inequalities with a constraint on the decision variables. 3) Constraint causes the magnitude of the eigenvalues to increase and aligns the impulse $\Delta v$ along the local stable manifold, a locally-optimal solution.}
    \label{fig:intro}
    \vspace{-0.8cm}
\end{figure}

Figure \ref{fig:intro} shows a visualisation of our method. First, a UPO is selected which will be stabilised, intersecting a local surface of section at the point $\mathbf{\Bar{x}}$. We then utilise our novel state-dependent augmentation strategy which enables the learning of state-dependent local Poincar\'e maps for volume-preserving systems. This simply consists of adding small velocity perturbations to the initial conditions of UPOs surrounding $\mathbf{\Bar{x}}$. These trajectories are simulated and successive crossings of a surface of section, $\boldsymbol{\Sigma}$, are recorded. From this data we discover a non-linear analytical expression of a Poincar\'e map using the SINDy algorithm, which we linearise around $\mathbf{\Bar{x}}$. In this work we also investigate where surfaces of section can be best placed, analysing control effectiveness and variance in the model coefficients using ensemble methods.
We stabilise the UPO by using a series of small velocity impulses, found by formulating the control problem as a series of linear matrix inequalities \cite{Parrilo2012SemidefiniteOptimization}. This is solved using an interior point method \cite{Nesterov1994Interior-PointProgramming}, where we add a constraint on the norm of the decision variables,  resulting in velocity impulses acting along the local stable manifold, as shown in Figure \ref{fig:intro}, significantly reducing the control cost.

Our method overcomes two key challenges. First, since conservative systems preserve phase-space volume, and therefore lack an underlying attractor, randomly initialised trajectories will not evolve to fall onto one common shape. Second, the original OGY method and subsequent data driven method in \cite{Bramburger2020Data-DrivenOrbits} relies on being able to change the system parameter. The augmentation and control strategy presented in our work addresses these challenges and enables quick-to-implement state-dependent stabilisation of UPOs of the 3BP. Critically, our method is sample efficient and discovers analytical expressions for local Poincar\'e maps, which are easily verified, a key advantage over black-box deep learning methods. In doing so, this facilitates the use of well developed linear control strategies. While our work herein is applied to celestial mechanics, future applications may include both the parameter and state-dependent control of other chaotic volume-preserving and dissipative systems, such as plasma, fluid flows, and robotics respectively.

\textbf{Outline.} In section \ref{sec: Chaos} we describe how our novel method discovers Poincar\'e maps and leverages linear matrix inequalities and a constrained projective method to find a controller that stabilises UPOs.
Following this, section \ref{sec: 3BP} describes the theory of the 3BP, how the stability of UPOs can be studied using the so called monodromy matrix, and the orbit families we focus on. In section \ref{sec: Results} we present our results and analysis of the discovery of local Poincar\'e maps and the stabilisation of UPOs of the 3BP. Finally, a discussion and summary of the key results is presented in section \ref{sec: Discussion}, as well as the limitations and avenues for future work.

\section{Chaotic Systems and their Control}

\label{sec: Chaos}
In this section we first discuss chaotic systems and their UPOs generally, and then detail our novel method for mapping discovery and control. We highlight how chaotic systems can be studied using Poincar\'{e} maps in \S \ref{sec: UPOs} and in \S \ref{sec: Learning maps} we outline our method for mapping discovery, a new augmentation strategy in conjunction with the SINDy algorithm, which enables the discovery of analytical expressions of these maps. Finally, in \S \ref{sec: Control theory} we detail how the pole-placement control problem is formulated by utilising a series of linear matrix inequalities, how a constraint was added, and an interior point method used to solve this convex problem.

\subsection{Unstable Periodic Orbits}
\label{sec: UPOs}
Chaotic systems are characterised by three key properties. First, their sensitivity  to small perturbations, where two trajectories starting close together will diverge exponentially with time. Second, they contain infinitely many recurrent trajectories known as UPOs. Third, they typically contain ergodic trajectories which, if left to evolve, will eventually visit a small neighbourhood of any UPO embedded in the attractor \cite{Ott2002ChaosSystems}. This final property is not always so apparent in conservative systems. A key distinction between dissipative and conservative systems is that the latter are phase-space volume-preserving, i.e. for a given set of initial conditions left to evolve, the volume that they occupy will remain the same. As such, the absence of a regular attractor structure in conservative systems necessitates the development of a novel augmentation method to discover accurate maps, which is presented.

The essence of a so called Poincar\'e map is to reduce a continuous-time system to a lower-dimensional discrete mapping. Considering the flow, $\mathbf{x}(t)$, of a continuous-time system, a lower dimensional surface of section, $\Sigma$, known as a Poincar\'{e} section can be defined transverse to the flow, with the flow crossing in only one direction. 
Successive crossings of trajectories with $\Sigma$ are considered, defining the Poincar\'e map
\begin{equation}
    \mathbf{x}_{n+1} = \mathbf{F}(\mathbf{x}_n),
    \label{eq: OG map}
\end{equation}
where $\mathbf{x}_n \in \Sigma$ at the $n$th intersection, and $\mathbf{x}_{n+1} \in \Sigma$ at the next. Periodic orbits of the original continuous time system occur as fixed points or cyclic points of the mapping itself. Often $\Sigma$ can be placed arbitrarily, so long as it is transverse to the flow, and there is no well-known rule as to where is the best location, at least for control purposes.

Throughout this work, we refer to the stability of an orbit and we define this based on the following conditions. Suppose there exists a fixed point of the map, $\Bar{\mathbf{x}}$, such that $\Bar{\mathbf{x}}=\mathbf{F}(\Bar{\mathbf{x}})$. Consider a small deviation from this point, $\boldsymbol{\delta}_n$, such that $\mathbf{x}_n = \Bar{\mathbf{x}} + \boldsymbol{\delta}_n$, and a Taylor expansion such that
\begin{equation}
    \mathbf{x}_{n+1} = \Bar{\mathbf{x}} + \boldsymbol{\delta}_{n+1} = \mathbf{F}(\Bar{\mathbf{x}}) + \mathbf{F}_x(\Bar{\mathbf{x}}) \boldsymbol{\delta}_n + \text{higher order terms}.
\end{equation}
where $\mathbf{F}_x(\Bar{\mathbf{x}})$ is the Jacobian of $\mathbf{F(x)}$ evaluated at $\Bar{\mathbf{x}}$. Since $\Bar{\mathbf{x}}$ is a fixed point this gives the result
\begin{equation}
    \boldsymbol{\delta}_{n+1} \approx  \mathbf{F}_x(\Bar{\mathbf{x}}) \boldsymbol{\delta}_n,
    \label{eq: xn+1}
\end{equation}
which is equivalent to 
\begin{equation}
    \mathbf{x}_{n+1} - \Bar{\mathbf{x}} \approx \mathbf{A}(\mathbf{x}_{n} - \Bar{\mathbf{x}}),
    \label{eq: lin sys}
\end{equation}
where $\mathbf{A}=\mathbf{F}_x (\Bar{\mathbf{x}})$. 
The linearised map, $\mathbf{A}$, describes how the deviation from $\Bar{\mathbf{x}}$ grows from one crossing of the section to the next, and therefore can be used to describe the stability of the periodic orbit. As such, we consider a matrix $\mathbf{M}$, and define
\begin{equation}
   \mathbf{M} \text{ is stable if } \exists \mathbf{P} \succ 0 : \mathbf{P} - \mathbf{M}^\top \mathbf{PM} \succ 0,
   \label{eq: Lyap stab}
\end{equation}
where superscript $\top$ is the matrix transpose and $\mathbf{P} \succ 0$ denotes that $\mathbf{P}$ is positive definite. This defines the notion of asymptotic stability, which can equally be stated in terms of the eigenvalues of a matrix: $\mathbf{M}$ is stable if all its eigenvalues $\lambda$ are such that $|\lambda| < 1$.

\subsection{Learning Poincar\'{e} Maps for Conservative Systems}
\label{sec: Learning maps}
In all but the simplest cases it is not possible to express a Poincar\'e map explicitly and so data-driven methods are required to obtain an analytical expression. Here we present a novel and straightforward strategy used in conjunction with the SINDy algorithm to discover such mappings, where we utilise prior knowledge of UPOs and augment their states. This strategy was developed due to the volume-preserving nature of conservative dynamical systems and the resulting lack of attractor structure, which necessitated a new approach for seeding the mapping discovery.

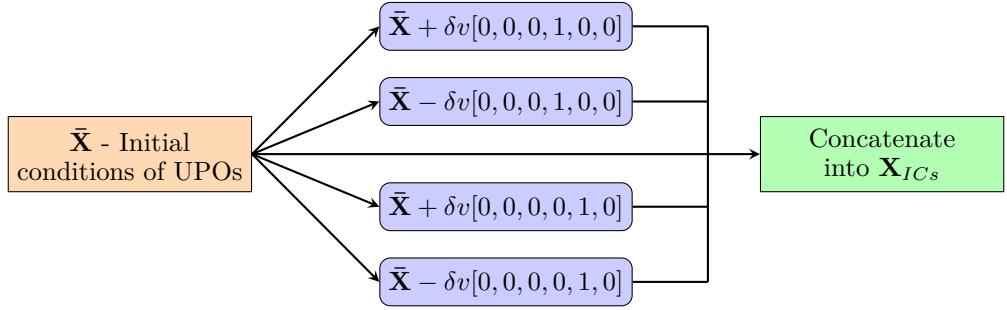
\begin{figure}[h!]
    \centering
    
    \begin{tikzpicture}[node distance=2cm]
    
    \node (start) [init, text width=3cm] {$\mathbf{\Bar{X}}$ - Initial conditions of UPOs};
    
    \node (n1a) [small, right of = start, xshift=3cm, yshift=1.7cm] {$\mathbf{\Bar{X}} + \delta v [0,0,0,1,0,0]$};
    
    \node (n1b) [small, below of = n1a, yshift=1.0cm] {$\mathbf{\Bar{X}} - \delta v [0,0,0,1,0,0]$};
    
    \node (n1c) [small, below of = n1b, yshift=0.6cm] {$\mathbf{\Bar{X}} + \delta v [0,0,0,0,1,0]$};
    
    \node (n1d) [small, below of = n1c, yshift=1.0cm] {$\mathbf{\Bar{X}} - \delta v [0,0,0,0,1,0]$};
    
    \node (n2) [decision, right of = start, xshift = 8cm, text width=3cm] {Concatenate into $\mathbf{X}_{ICs}$};
    
    \draw [arrow] (start.east) to (n1a.west);
    \draw [arrow] (start.east) to (n1b.west);
    \draw [arrow] (start.east) to (n1c.west);
    \draw [arrow] (start.east) to (n1d.west);
    
    \coordinate (m1) at ($(n1a.east)+(1,0)$);
    \coordinate (m2) at ($(n1b.east)+(1,0)$);
    \coordinate (m3) at ($(n1c.east)+(1,0)$);
    \coordinate (m4) at ($(n1d.east)+(1,0)$);
    
    \draw [line] (n1a.east) -- (m1);
    \draw [line] (n1b.east) -- (m2);
    \draw [line] (n1c.east) -- (m3);
    \draw [line] (n1d.east) -- (m4);
    
    \draw [line] (m1) -- (m4);
    
    \coordinate (mMid) at ($(m1)!0.5!(m4)$);
    
    \draw [arrow] (mMid) -- (n2.west);

    \draw[line] (start.east) -- (mMid);

    \end{tikzpicture}
    \caption{\textbf{Augmentation Strategy -} We take the initial conditions of the desired UPO and $m$ UPOs centred around it in energy level, concatenating them into $\mathbf{\Bar{X}}$. Next, we add small positive and negative velocity perturbations $\delta v{}$ in $x$ and $y$ separately, to each UPO. These augmented initial conditions and the original $\mathbf{\Bar{X}}$ are concatenated into one initial conditions vector $\mathbf{X}_{ICs}$.}   
    \label{fig: aug diag}
\end{figure}

Considering the flow of the system, we define a surface of section, $\Sigma$, as described in \S \ref{sec: UPOs}. We take the initial state vector of the target UPO, and we consider small perturbations in energy level (see \S \ref{sec: 3BP eqns} for further detail) above and below to find $m$ UPOs nearby, which are recorded in a vector, denoted $\mathbf{\Bar{X}}$ in Figure \ref{fig: aug diag}. Small velocity perturbations, $\delta v$, are applied to these initial conditions and combined into one matrix, as shown in Figure \ref{fig: aug diag}. The size of the velocity perturbation is not known \textit{a priori}, and this is adjusted in order to discover an accurate mapping.
The mapping discovery is seeded by simulating these $j= \{1,...,5m\}$ different initial conditions, with an 8th order Runge-Kutta method \cite{Govorukhin2024Ode87Integrator}. Successive crossings of $\Sigma$ are recorded, yielding the sequence $\{\mathbf{x}_n = \mathbf{x}(t_n) : \mathbf{x}(t_n) \in \Sigma, t_n > t_{n-1}\}$ which results in the data series $\{ \mathbf{x}^{(j)}_1,..., \mathbf{x}^{(j)}_{N_{j}}\}$  where $N_j$ is the number of crossings of $\Sigma$ \cite{Bramburger2020Data-DrivenOrbits}. Since we are performing a local mapping discovery, only crossings within a local region which satisfy
\begin{equation}
    || \mathbf{x}_n - \Bar{\mathbf{x}}|| \leq \eta
    \label{eq: eta}
\end{equation}
are selected, where $\eta$ is a user defined threshold value. Two vectors are created
\begin{equation}
\mathbf{X}_1 =
\begin{bmatrix}
  \mathbf{x}^{(1)}_1 \\
  \vdots \\
  \mathbf{x}^{(1)}_{N_1-1} \\
  \mathbf{x}^{(2)}_1 \\
  \vdots \\
  \mathbf{x}^{(2)}_{N_2-1} \\
  \vdots \\
  \mathbf{x}^{(M)}_1 \\
  \vdots \\
  \mathbf{x}^{(M)}_{N_M-1}
\end{bmatrix},
\qquad
\mathbf{X}_2 =
\begin{bmatrix}
  \mathbf{x}^{(1)}_2 \\
  \vdots \\
  \mathbf{x}^{(1)}_{N_1} \\
  \mathbf{x}^{(2)}_2 \\
  \vdots \\
  \mathbf{x}^{(2)}_{N_2} \\
  \vdots \\
  \mathbf{x}^{(M)}_2 \\
  \vdots \\
  \mathbf{x}^{(M)}_{N_M}
\end{bmatrix},
\end{equation}
where for each entry in $\mathbf{X}_1$ the corresponding one in $\mathbf{X}_2$ is the next crossing of the section. 
The goal is to discover a mapping, $\mathbf{F}$, from this data, which maps each row of $\mathbf{X}_1$ onto the same row in $\mathbf{X_2}$, following the same form as Equation \eqref{eq: OG map}. We select a library of candidate functions, $\{\theta_1,...,\theta_D\}$, and seek to identify the mapping within the linear span of the library, given by
\begin{equation}
    \mathbf{F} = \xi_1 \theta_1(\mathbf{x}) + ... + \xi_D \theta_D(\mathbf{x}),
    \label{eq: lin span}
\end{equation}
where, $\boldsymbol{\Xi} = \xi_1,...,\xi_D$, are the coefficients of the linear span. Combining Equations \eqref{eq: OG map} and \eqref{eq: lin span} gives
\begin{equation}
    \mathbf{x}_{n+1} = \xi_1 \theta_1(\mathbf{x}_n) + ... + \xi_D \theta_D(\mathbf{x}_n),
\end{equation}
which can also be expressed as 
\begin{equation}
    \mathbf{X}_2= \boldsymbol{\Theta}(\mathbf{X}_1)\boldsymbol{\Xi},
\end{equation}
where $\boldsymbol{\Theta}$ is the library of candidate functions evaluated at each entry of $\mathbf{X}_1$.
The SINDy algorithm utilises sparse regression methods to find $\boldsymbol{\Xi}$, since many physical systems can be expressed with a small number of candidate functions, which helps prevent overfitting to measurement error. In this paper we use the sequentially thresholded least squares regression (STLSR) optimiser due to its sharp rate of convergence \cite{Zhang2019OnAlgorithm}. More specifically, the following iterative procedure is used: 
\begin{enumerate}
    \item Least squares regression:  \begin{equation}
    {\boldsymbol{\Xi}} = \argmin_{\Xi}||\mathbf{X}_2 - \boldsymbol{\Theta}(\mathbf{X}_1)\boldsymbol{\Xi} ||_2
    \end{equation} 
    \item Thresholding: each $\xi_j$ term in $\boldsymbol{\Xi}$ that satisfies  $|\xi_j| < \lambda_{\text{sparse}}$ is set to zero, where $\lambda_{\text{sparse}} > 0$ is a user-defined sparsity promoting parameter. 
\end{enumerate}
These two steps are then repeated until all $\xi_j \geq \lambda_{\text{sparse}}$ and a sufficiently sparse model is obtained. Selecting the candidate library functions and $\lambda_{\text{sparse}}$ are key to obtaining an accurate model, and the two are related since the closer the candidate library to the true terms, the smaller $\lambda_{sparse}$ can be. In the case of Poincar\'{e} maps we do not know the true map, so a simple order five polynomial library is used \cite{Bramburger2020Data-DrivenOrbits}, which proved to be effective.

Throughout this paper, mappings are always centred and linearised around the UPO which we wish to stabilise. Considering the desired crossing, $\Bar{\mathbf{x}}$, of the UPO with $\Sigma$, the state error vector is defined, $\mathbf{\Tilde{x}}_n = \mathbf{x}_n - \mathbf{\Bar{x}}$, which is the same as $\boldsymbol{\delta}_n$ in Equation \eqref{eq: xn+1}. The discovered mapping is therefore of the form 
\begin{equation}
    \mathbf{\Tilde{x}}_{n+1} = \mathbf{\Tilde{F}}(\mathbf{\Tilde{x}}_n). 
    \label{eq: tilde sys}
\end{equation}
This does not change the preceding stability analysis of the Poincar\'{e} map nor the control formulation that follows, since this is based on the linearised system expressing how the error grows. That is to say, in the linearised system \eqref{eq: lin sys}, $\mathbf{A}=\mathbf{\Tilde{F}}_{\Tilde{x}}(\mathbf{\mathbf{0}})$ as well, since the zero vector is a fixed point of the system \eqref{eq: tilde sys}.

\textbf{Data-ensembling.} We also employ ensemble learning via eSINDy, which provides model-variance estimates and is more robust to noise \cite{Fasel2021Ensemble-SINDy:Control}. Specifically, we perform data bootstrapping where we take random samples of the data and apply standard SINDy to each, discovering $N_e$ candidate models. Coefficients with an inclusion probability less than a user-defined threshold $\rho$ are thresholded out, and the median model is taken as the final model. Importantly, calculating the variance $\sigma_i$ of each coefficient, $\xi_i$, gives us insight into the model uncertainties. For a full explanation of the method we direct the reader to \cite{Fasel2021Ensemble-SINDy:Control}.

\subsection{State-Dependent Stabilisation of UPOs}
\label{sec: Control theory}
The original OGY method relies on parameter perturbations to stabilise UPOs \cite{Ott1990ControllingChaos}; however, since  many real world systems allow control of only specific states rather than parameters, this paper focuses on the state-dependent scenario. Here, an outline of a state-dependent formulation of the problem formulation is provided.
Small impulses are applied to the system at each crossing of $\Sigma$, based on the learned Poincar\'e map
\begin{equation}
    \mathbf{x}_{n+1} = \mathbf{F}(\mathbf{x}_n + \mathbf{Ru}_n),
    \label{eq: control map}
\end{equation}
where $\mathbf{R}$ represents the physical restrictions of the control inputs. The aim is to find small control inputs $\mathbf{u}_n$ to stabilise a fixed point $\Bar{\mathbf{x}}$, of the mapping. We consider state-dependent feedback control that helps to re-orient the trajectory back towards the UPO according to
\begin{equation}
    \mathbf{u}_n = 
    \begin{cases}
    \mathbf{K}(\mathbf{x}_n - \Bar{\mathbf{x}}) & || \mathbf{x}_n - \Bar{\mathbf{x}}|| \leq \eta \\
    \mathbf{0} & || \mathbf{x}_n - \Bar{\mathbf{x}}|| > \eta \\
    \end{cases},
    \label{eq: cont}
\end{equation}
where $\eta$ is a user-defined threshold parameter, chosen such that $|\mathbf{K}(\mathbf{x}_{n} - \Bar{\mathbf{x}} )|$ is sufficiently small, guaranteeing the trajectories remain within the region of validity for the linearised system $\mathbf{A}$. We linearise the mapping \eqref{eq: control map} around $\Bar{\mathbf{x}}$ and $\mathbf{u}_n = 0$, resulting in the equation
\begin{equation}
    \mathbf{x}_{n+1} - \Bar{\mathbf{x}} \approx \mathbf{A}(\mathbf{x}_{n} - \Bar{\mathbf{x}}) + \mathbf{ARK}(\mathbf{x}_{n} - \Bar{\mathbf{x}} ) = (\mathbf{A} + \mathbf{ARK}) (\mathbf{x}_{n} - \Bar{\mathbf{x}} ).
    \label{eq: A+BK}
\end{equation}

The problem is now to find a matrix, $\mathbf{K}$, such that the orbits are driven to or remain near the desired UPO. In \S \ref{sec: UPOs} we described the criterion for a matrix to be considered Lyapunov stable. Taking the definition from \eqref{eq: Lyap stab} for the linearised controlled system, we consider $\mathbf{M = A+BK}$ where $\mathbf{B = AR}$, which gives the condition
\begin{equation}
    \mathbf{P} - (\mathbf{A+BK})^\top \mathbf{P(A+BK)} \succ 0,
    \label{eq: Lyap Eqn}
\end{equation}
which is nonlinear with respect to the unknown $\mathbf{K}$. More generally, finding a $\mathbf{K}$ such that the system is stable is known as the pole-placement problem, and is only feasible if the controllability matrix
\begin{equation}
    \mathbf{C} = (\mathbf{B} : \mathbf{A}\mathbf{B} : \mathbf{A}^2\mathbf{B} : \dots : \mathbf{A}^{n-1}\mathbf{B})
\end{equation}
is full rank.

\subsubsection{Convex Optimisation Formulation}
To determine the controller matrix $\mathbf{K}$ we turn to the method described in \cite{Parrilo2012SemidefiniteOptimization} which reformulates \eqref{eq: Lyap Eqn} as a series of linear matrix inequalities (LMIs). The advantage is that LMIs are convex semidefinite optimisation problems which can be solved quickly and efficiently, by using interior point methods implemented within a variety of software packages. Taking the Schur complement of Equation \eqref{eq: Lyap Eqn} we obtain the equivalent expression 
\begin{equation}
    \begin{bmatrix}
    \mathbf{P} & (\mathbf{A}+\mathbf{B}\mathbf{K})^\top \mathbf{P} \\
    \mathbf{P}(\mathbf{A}+\mathbf{B}\mathbf{K}) & \mathbf{P}
\end{bmatrix} \succ 0,
\label{eq: Schur}
\end{equation}
which is not yet convex as it is bilinear in $(\mathbf{K,P})$. Recalling that $\mathbf{P}$ is assumed to be strictly positive definite, it is invertible and so we define $\mathbf{Q = P}^{-1}$. Multiplying \eqref{eq: Schur} on the left and right by the block diagonal $\text{diag}(\mathbf{Q,Q})$ gives
\begin{equation}
    \begin{bmatrix}
    \mathbf{Q} & \mathbf{Q}(\mathbf{A}+\mathbf{B}\mathbf{K})^\top  \\
    (\mathbf{A}+\mathbf{B}\mathbf{K})\mathbf{Q} & \mathbf{Q}
\end{bmatrix} \succ 0.
\notag
\end{equation}
Again, defining a new variable $\mathbf{Y=KQ}$ results in the problem
\begin{equation}
    \begin{bmatrix}
    \mathbf{Q} & \mathbf{QA}^\top + \mathbf{Y}^\top\mathbf{B}^\top \\
    \mathbf{AQ + BY} & \mathbf{Q}
\end{bmatrix} \succ 0,
\label{eq: LMI form}
\end{equation}
which is linear in the variables $(\mathbf{Q,Y})$ and forms an auxiliary convex programming problem. The solution $\mathbf{K}$ exists so long as the set of the decision variables $(\mathbf{Q,Y})$ solving \eqref{eq: LMI form} is non-empty, and is recovered by
\begin{equation}
    \mathbf{K = Y Q}^{-1}.
    \label{eq: K sol}
\end{equation}
To implement this method numerically, we used the \lstinline{feasp} optimiser in MATLAB based on Nesterov and Nemirovski's Projective Method which converges with polynomial time \cite{Nesterov1994Interior-PointProgramming}.
Critically, in our method we constrain the magnitude of the decision variables by defining a feasibility radius $R$. The sum of the squares of the decision variables are constrained such that
\begin{equation}
    \sum_{i=1}^{N} d_i^2 < R^2, 
    \label{eq: R}
\end{equation}
where $d$ contains the elements of $\mathbf{Q}$ and $\mathbf{Y}$, and so is equivalent to $||\mathbf{Q}||_2 + ||\mathbf{Y}||_2 < R$. This constraint reduces the size of the feasible set but $R$ can only be reduced so long as the problem \eqref{eq: LMI form} remains feasible.

\section{The Three-Body Problem}
\label{sec: 3BP}
The three-body problem (3BP) is a conservative chaotic system which describes the movement of three bodies that exert attractive forces on each other. In this study we analyse the Earth-Moon 3BP, where the attractive force is gravity, but the 3BP appears in many different scientific fields \cite{Dzhanoev2010ChaosProblem}. The 3BP was first studied in depth by Henry Poincar\'e in the 19th century, who proved periodic orbits exist within the 3BP, and is still the subject of much ongoing research. 
This section outlines the governing equations and dynamics of the 3BP. We briefly explain invariant manifold theory and how it is used to categorise the stability of a UPO as well as the local manifold structure. Finally, we discuss the two orbit families that we use in this paper and where the different surfaces of section are placed. 

\subsection{Equations of Motion}
\label{sec: 3BP eqns}
In this paper we study the circular restricted 3BP, abbreviated CR3BP. This describes the motion of a negligible-mass particle, $P_3$, under the influence of gravitational fields of two larger masses, $m_1$ and $m_2$, which follow circular orbits around their common centre of mass. In celestial mechanics this particle is often a comet or spacecraft.
The mass ratio 
\begin{equation}
    \mu = \frac{m_2}  {m_1 + m_2},
\end{equation}
is the only system parameter, where $m$ are the masses and $m_1 > m_2$. Similar to Reynolds number in fluid mechanics, this determines the types of behaviour observed. In this paper we consider the Earth-Moon CR3BP which has a mass ratio $\mu = 1.215059\times10^{-2}$.

\begin{figure}[htb]
    \includegraphics[width=9cm]{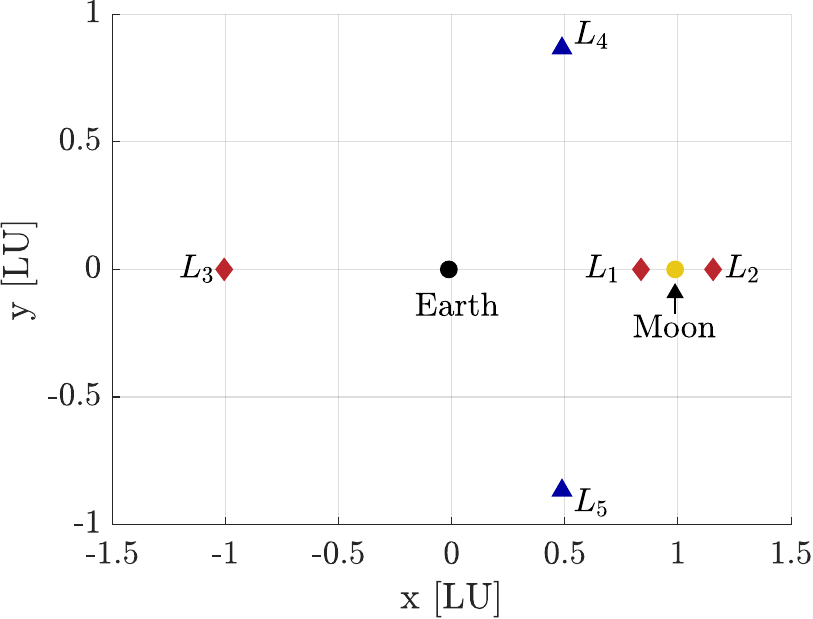}
     \centering
    \captionof{figure}{Location of the unstable $L_1, L_2, L_3$ and stable $L_4,L_5$ Lagrange points of the Earth-Moon CR3BP in the orbital plane.}
    \label{fig: Lagrange}
\end{figure}

The common centre of mass of the two larger bodies defines the origin of a rotating reference frame, which is used throughout as orbit patterns are shown more intuitively \cite{Koon2011DynamicalDesign.}. In this coordinate system, the bodies $m_1$ and $m_2$ are located at $(-\mu,0,0)$ and $(1-\mu,0,0)$ respectively and the $x$-axis points towards $m_2$, as shown in Figure \ref{fig: Lagrange}. The larger bodies orbit in only the $x$-$y$ orbital plane and $x$-$y$ axes rotate with respect to the inertial frame at the same angular rate as the body $m_2$. The motion of $P_3$ is governed by the equations
\begin{equation} 
    \begin{aligned} 
    \ddot{x} - 2\dot{y} - x &= - \frac{(1 - \mu)(x + \mu)}{d^3} - \frac{\mu}{r^3} (x - 1 + \mu) \\ 
    \ddot{y} + 2\dot{x} - y &= - \frac{(1 - \mu)}{d^3} y - \frac{\mu}{r^3} y \\ 
    \ddot{z} &= - \frac{(1 - \mu)}{d^3} z - \frac{\mu}{r^3} z 
    \end{aligned} ,
    \label{eq:CR3BP} 
\end{equation}
where $d$ and $r$ are the distance from $P_3$ to $m_1$ and $m_2$ respectively, calculated as
\begin{equation}
    \begin{aligned}
        d &= \sqrt{ (x+\mu)^2 + y^2 + z^2} \\
    r &= \sqrt{ (x-1 +\mu)^2 + y^2 + z^2}
    \end{aligned}.
\end{equation}
These equations are given using non-dimensional units, which are used throughout this paper. A length unit, LU, is the distance between $m_1$ and $m_2$, of 389703 km, and the time unit, TU, is the inverse of the relative angular frequency between the two large bodies, of 382981 seconds.

The CR3BP admits a conserved energy integral known as the Jacobi constant,
\begin{equation}
    C = 2U - V^2,
    \label{eq: C}
\end{equation}
where $V = \sqrt{\dot{x}^2 + \dot{y}^2 + \dot{z}^2}$, and $U$ is the effective potential experienced by $P_3$ given by
\begin{equation}
    U = \frac{x^2 + y^2}{2} + \frac{\mu}{r} + \frac{1-\mu}{d}.
    \label{eq: pot}
\end{equation}
The Jacobi constant describes the total energy of the particle $P_3$, which constrains the possible trajectories of a particle to a certain region of space known as a Hill's region. Since this is the only integral constraining the particle's motion, the problem is therefore non-integrable. Due to the conservative nature of the 3BP we use an 8th order Runge-Kutta integrator \cite{Govorukhin2024Ode87Integrator} with a relative tolerance of $10^{-12}$ and absolute tolerance $10^{-14}$, to avoid large deviations in the Jacobi constant.

The CR3BP admits five equilibria where the gravitational pull of $m_2$ and $m_1$ are in balance with the fictitious accelerations. They are known as Lagrange or libration points, and for the value of $\mu$ considered, their locations in the $x$-$y$ plane are given in Table \ref{tab: Lagrange}. The $z$ value of each Lagrange point is zero and therefore not given in Table \ref{tab: Lagrange} for brevity.
\captionof{table}{Location and Jacobi constant of Lagrange points in the Earth-Moon system.}
\centering
\label{tab: Lagrange}
\begin{tabular}{cccccc}
\toprule
 & $L_1$ & $L_2$ & $L_3$ & $L_4$ & $L_5$ \\
\midrule
x (LU) & 0.837 & 1.156 & -1.005 & 0.488 & 0.488 \\
\midrule
y (LU) & 0.000 & 0.000 & 0.000 & 0.866 & -0.866 \\
\midrule
C (LU$^2$/TU$^2$) & 3.188 & 3.172 & 3.013 & 2.988 & 2.988 \\
\bottomrule
\end{tabular}
\justify
The locations of the Lagrange points in the orbital plane of the Earth-Moon CR3BP are presented in Figure \ref{fig: Lagrange}. The collinear Lagrange points $L_1, L_2, L_3$, are unstable saddle points whereas $L_4$ and $L_5$ are stable due to the effect of the Coriolis force \cite{Koon2011DynamicalDesign.}. 

Various families of both stable and unstable orbits exist; some orbit the Lagrange points or jump between them, whilst others, known as resonant orbits, have a period commensurate to that of the system. Poincar\'{e} proved that UPOs exist within the CR3BP and modern computational methods have enabled the discovery of a plethora of different UPO families \cite{Baresi2020AProblems, Doedel2011ComputationProblem, Dzhanoev2010ChaosProblem, Restrepo2018ASystem}. In this work, we consider UPOs that orbit the Lagrange points, where an orbit family consists of infinitely many UPOs with values of $C$ slightly below that of the associated Lagrange point.

\subsection{Invariant Manifold Theory}
\label{sec: Invariant M}
Throughout this paper, we have referred to the stability of periodic orbits, namely UPOs. The goal of this section is to provide the relevant theory for calculating the stability properties of a UPO, in the context of the CR3BP. Consider a dynamical system
\begin{equation}
    \mathbf{\dot{x} = f(x)}, \quad \mathbf{x} \in \mathbb{R}^p
    \label{eq: basic sys}
\end{equation}
with a flow map $\boldsymbol{\phi}(t,\mathbf{x}_0) : \mathbf{x}(t_0) \rightarrow \mathbf{x}(t)$, which maps the state of a particle from time $t_0$ to $t$. Considering a perturbed initial condition $\Bar{\mathbf{x}}_0 + \delta\Bar{\mathbf{x}}_0$, the displacement from a reference trajectory ${\mathbf{\Bar{x}}}(t)$, is given by
\begin{equation}
    \delta\Bar{\mathbf{x}}(t) = \boldsymbol{\phi}(t; \Bar{\mathbf{x}}_0 + \delta\Bar{\mathbf{x}}_0) - \boldsymbol{\phi}(t;\Bar{\mathbf{x}}_0),
    \notag
\end{equation}
which, when $\delta$ is small, we can expand using the Taylor series
\begin{equation}
     \delta\Bar{\mathbf{x}}(t) = \frac{\partial \boldsymbol{\phi}(t;\Bar{\mathbf{x}}_0)}{\partial \Bar{\mathbf{x}}_0} \delta \mathbf{\Bar{x}}_0 + \text{higher order terms}.
     \notag
\end{equation}
The term $\frac{\partial \boldsymbol{\phi}(t;\Bar{\mathbf{x}}_0)}{\partial \Bar{\mathbf{x}}_0}$ is known as the state transition matrix, denoted $\boldsymbol{\Phi}(t,t_0)$, which describes the linearised relationship between small initial perturbations at $t_0$ and the corresponding change in state at the final time $t$. In order to compute the state transition matrix, the initial value problem
\begin{equation}
    \Dot{\boldsymbol{\Phi}}(t,t_0) = \mathbf{f_x}(\Bar{\mathbf{x}}(t)) \boldsymbol{\Phi}(t,t_0), \quad \text{where } \boldsymbol{\Phi}(t_0,t_0) = \mathbf{I},
    \label{eq: STM sys}
\end{equation}
is considered, where $\mathbf{f_x}$ denotes the Jacobian of $\mathbf{f}$ and $\mathbf{I}$ is the identity matrix. In all but the simplest of scenarios, this problem cannot be solved analytically since since we do not have an explicit analytical description of the solution ${\mathbf{\Bar{x}}}(t)$.
We must therefore numerically integrate Equation \eqref{eq: STM sys}, which gives a system of $p^2$ first-order scalar differential equations expressing the elements of $\boldsymbol{\Phi}$. In the case of the CR3BP, this gives a system of 36 differential equations, since there are six states $(x, y, z, \dot{x}, \dot{y}, \dot{z})$. Combining Equations \eqref{eq: basic sys} and \eqref{eq: STM sys} results in a system of 42 differential equations which can be integrated concurrently to calculate the state transition matrix.

Let us now consider the particular case of the reference trajectory being a periodic orbit, i.e. $\mathbf{x}(t + T) = \mathbf{x}(t)$ for a given $T > 0$ and $t \in \mathbb{R}$. The resulting displacement after one period, $T$, is therefore given by
\begin{equation}
    \delta\Bar{\mathbf{x}}(T) = \boldsymbol{\Phi}(T)\delta\Bar{\mathbf{x}}_0,
\end{equation}
to first order accuracy. The state transition matrix after one period, $\boldsymbol{\Phi}(T)$, is referred to as the monodromy matrix,
\begin{equation}
    \mathbf{M} = \frac{\partial \boldsymbol{\phi}(T;\Bar{\boldsymbol{x}}_0)}{\partial \boldsymbol{x}_0}.
    \label{eq: Monodromy}
\end{equation}

The stability of a given trajectory can be analysed similarly to the method outlined in \S \ref{sec: Chaos} for the linearised Poincar\'{e} map, since $\mathbf{M}$ is also a linearised map. 
The eigenvalues $\lambda$ of the monodromy matrix, known as Floquet multipliers, describe the stability of the trajectory according to \eqref{eq: Lyap stab}. The associated eigenvectors describe the local invariant manifolds of the system, where the stable manifold $\nu_s$ is associated with $\lambda < 1$ and the unstable manifold $\nu_u$ with $\lambda > 1$.
Notably, for a Hamiltonian system, the monodromy matrix will contain the eigenvalue $\lambda=1$ with algebraic multiplicity of at least two \cite{Koon2011DynamicalDesign.}, and the eigenvalues of $\mathbf{M}$ should be the same as those of the linearised Poincar\'e map \cite{Koon2011DynamicalDesign.}. 

Recalling that the monodromy matrix $\mathbf{M}$ represents the linearised evolution of small perturbations over one period, we use the norm $||\mathbf{M}||$ as a measure of sensitivity in this paper. The monodromy matricies at different surfaces of section are conjugate, so have the same eigenvalues and differ only by a change of basis, but comparing their norms highlights which section is more sensitive.

\subsection{Orbit Families}
There are many different families of periodic orbits within the CR3BP and to this day new orbits of multi-body systems continue to be discovered, as highlighted in \S \ref{sec: 3BP eqns}. We direct the reader to \cite{Vaquero2018Poincare:Tool} for an overview of the various families and the initial conditions of these UPOs in the six-dimensional ODE \eqref{eq:CR3BP}. In this paper, we focus on two of the most commonly used families of UPOs in space missions, Lyapunov and halo orbits, and we provide a control strategy to stabilise these UPOs. 

\subsubsection{Lyapunov Orbits}
\label{sec: Lyap fam}
\begin{figure}[h!]
\centering
\begin{minipage}{0.48\textwidth}
  \centering
  \includegraphics[width=.9\linewidth]{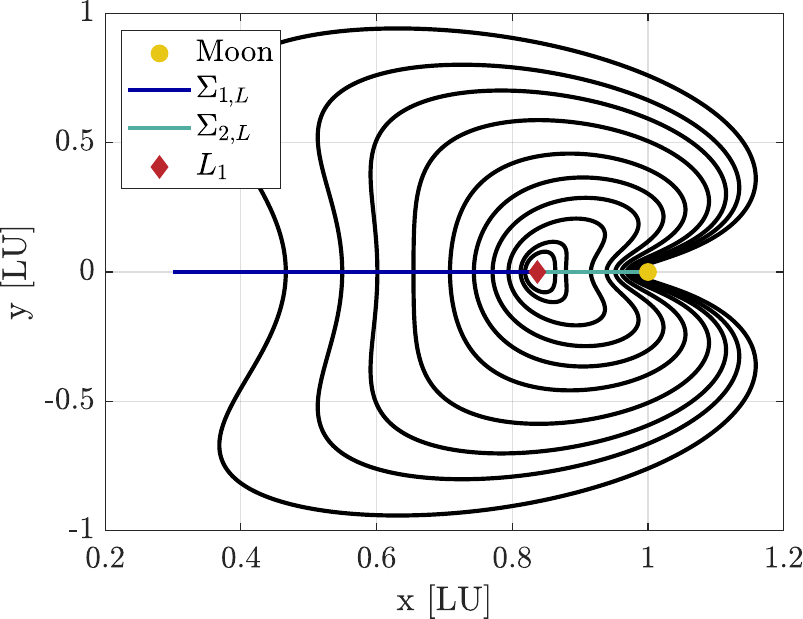}
  \captionof{figure}{Various $L_1$ Lyapunov orbits shown in the $x$-$y$ plane intersecting with the two surfaces of section, $\Sigma_{1,L}$ and $\Sigma_{2,L}$, we sample data from.}
  \label{fig:Lyap Fam}
\end{minipage}
\hfill
\begin{minipage}{.48\textwidth}
  \centering
  \includegraphics[width=.9\linewidth]{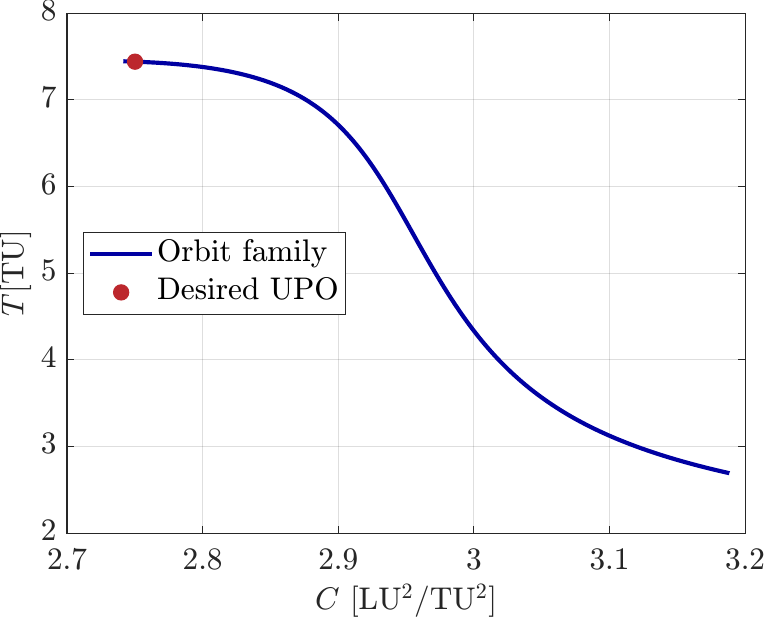}
  \captionof{figure}{Nonlinear inverse relationship between orbital period $T$ and Jacobi constant $C$ for the $L_1$ Lyapunov family, where all UPOs sit below the $L_1$ value of $C=3.1883$.}
  \label{fig:Lyap Props}
\end{minipage}
\end{figure}

The motion of Lyapunov orbits is restricted to the orbital plane of the two primary masses and as such are known as planar orbits. We focus on Lyapunov orbits around the $L_1$ Lagrange point, while similar methods can be applied to study those about $L_2$ as well. Figure \ref{fig:Lyap Fam} shows the orbit family with two different surfaces of section, defined as 
\begin{equation}
    \Sigma_{1,L} = \{ (x,\Dot{x},\Dot{y}) \text{ }|\text{ } x < 0.8369, z = \dot{z} = 0\},
    \label{eq:S1}
\end{equation}
\begin{equation}
    \Sigma_{2,L} = \{ (x,\Dot{x},\Dot{y}) \text{ }|\text{ } x > 0.8369, z = \dot{z} = 0\},
    \label{eq:S2}
\end{equation}
which satisfy the requirement that $\Sigma$ be transverse to the flow of the system, as discussed in \S \ref{sec: UPOs}. These surfaces have no $z$ component due to the motion of Lyapunov orbits occurring purely in the orbital plane. 
The $L_1$ Lagrange point has a Jacobi constant of $C=3.1883$ and the Lyapunov orbit family sits just below, spanning the range $C \in [2.7415, 3.1883]$.
Figure \ref{fig:Lyap Props} shows the nonlinear inverse relationship between the orbital period and Jacobi constant for the orbit family. A long-period orbit with $T = 7.4417$ [TU] and $C= 2.75018$ [TU$^2$/LU$^2$] was selected  as the desired UPO to be stabilised, shown within the orbit family in Figure \ref{fig:Lyap Props}. The monodromy matrix $\mathbf{M}$ of UPOs in the Lyapunov family, admit eigenvalues of the form
\begin{equation}
    \lambda_1 > 1, \quad \lambda_2 = \frac{1}{\lambda_1}, \quad \lambda_3 = \lambda_4 = 1,
    \label{eq: Lyap eigs}
\end{equation}
which occur as a reciprocal pair plus two unit values, due to the conservative (volume-preserving) nature of the system. The eigenvectors associated with $\lambda_1$ and $\lambda_2$ correspond to the local unstable and stable manifolds respectively. For the selected orbit the unstable and stable eigenvalues are $\lambda_u = 219.5457$ and $\lambda_s = 0.0046$ respectively.

\subsubsection{Halo Orbits}
\label{sec: Halo fam}
\begin{figure}[h!]
\centering
\begin{minipage}{0.48\textwidth}
    \centering
    \includegraphics[width=0.9\linewidth]{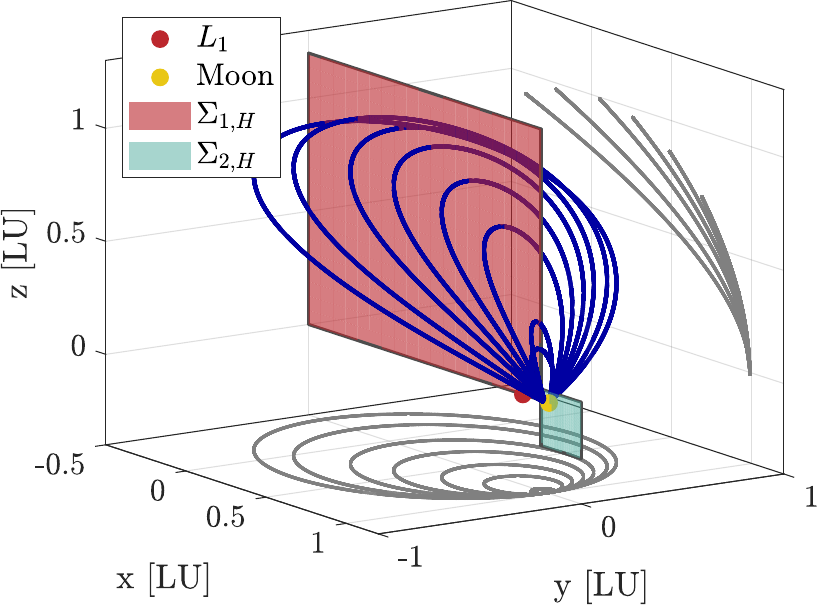}
    \caption{Various $L_1$ halo orbits intersecting transverse to the two surfaces of section, with projections onto the $x$-$y$ and $x$-$z$ planes showing highly concentrated region near the Moon, the larger primary body.}
    \label{fig:Halo Fam}
\end{minipage}
\hfill
\begin{minipage}{.48\textwidth}
  \centering
  \includegraphics[width=0.9\linewidth]{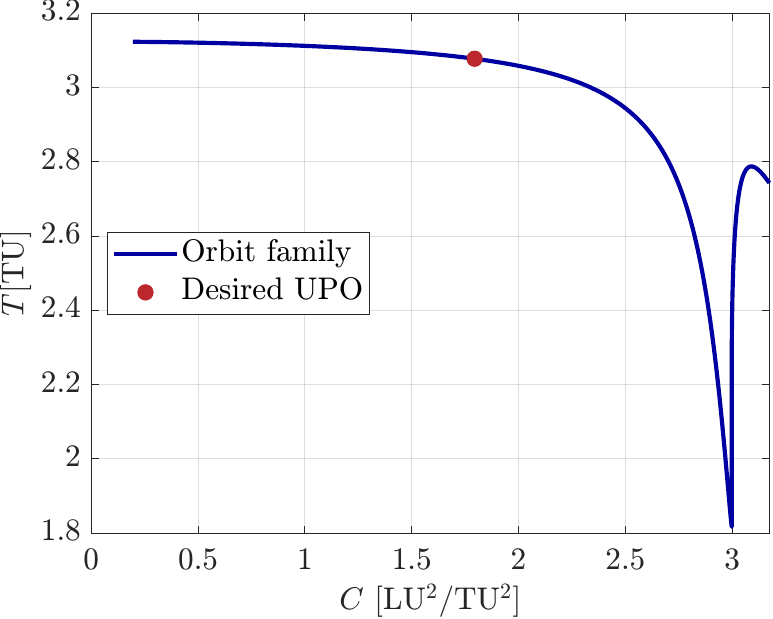}
    \caption{Highly nonlinear relationship between orbital period $T$ and Jacobi constant $C$ for the halo orbit family, where all UPOs sit below the $L_1$ value of $C=3.1883$. Desired UPO in shown in red.}
    \label{fig:Halo props}
\end{minipage}
\end{figure}

The halo orbit family consists of trajectories above and below the orbital plane, known as northern and southern orbits respectively. This paper focuses on $L_1$ northern halo orbits, which can be used to establish continuous links with the far side of the moon, but the methods presented can be applied to other orbits in a similar manner.
There are only two places where we can define $\Sigma$, satisfying the condition that they are transverse to the flow, given by:
\begin{equation}
    \Sigma_{1,H} = \{ (x,z,\Dot{x},\Dot{y},\Dot{z}) \text{ }|\text{ } z > 0, y=0 \},
    \label{eq:H1}
\end{equation}
\begin{equation}
    \Sigma_{2,H} = \{ (x,z,\Dot{x},\Dot{y},\Dot{z}) \text{ }|\text{ } z < 0, y=0 \}.
    \label{eq:H2}
\end{equation}
Various UPOs from the family are shown in Figure \ref{fig:Halo Fam}, with their projections onto the $x$-$y$ and $x$-$z$ planes, alongside the two surfaces of section.
Halo orbits span the range $C \in [0.1952, 3.1743]$ TU$^2$/LU$^2$ and their orbital periods generally remain within a smaller range than for Lyapunov orbits. There is a extreme change for $C > 2.5$, as seen in Figure \ref{fig:Halo props}, where the orbit shape changes markedly, but we do not focus on these in our paper. We selected the desired halo UPO with parameters $C=1.7979$ TU$^2$/LU$^2$ and $T=3.0773$ TU. The eigenvalues of the monodromy matrix $\mathbf{M}$ for halo orbits are of the form
\begin{equation}
    \lambda_1 > 1, \quad \lambda_2 = \frac{1}{\lambda_1}, \quad \lambda_3 = \lambda_4 = 1, \quad \lambda_5 = \Bar{\lambda}_6 \quad |\lambda_5| = 1,
    \label{eq: halo eigs}
\end{equation}
where $\lambda_1$ and $\lambda_2$ are the unstable and stable eigenvalues, which for the selected orbit are $\lambda_u = 480.2979$ and $\lambda_s = 0.0021$ respectively.

\section{Results}
\label{sec: Results}
In this section we apply our methods to successfully stabilise both Lyapunov and halo orbits, both UPOs of the 3BP. The augmentation strategy presented enables the discovery of accurate Poincar\'e maps of the three body problem. 
All the methods are implemented in MATLAB and the code for this is available at \url{https://github.com/ombrook/DDS_of_UPOs_of_3BP}. Throughout we use an order five polynomial library, the minimum degree for meaningful results. Discovering a linear map directly did not result in a map suitable for control, as this applies a global linear fit to all the data which does not capture the local dynamics around the target UPO accurately. Discovering a nonlinear map of sufficiently high order and then evaluating the Jacobian around the target UPO \eqref{eq: lin sys}, results in a truly local linearisation. Increasing the order beyond five did not cause any higher order terms to appear. Additionally, we found a sparsity promoting parameter $\lambda_{sparse} = 10^{-6}$ in the standard SINDy case, and a threshold parameter $\eta = 1.0$, lead to accurate model discovery. We invite the reader to investigate these properties for themselves using the accompanying code. The catalogue of initial conditions of all the UPOs can be found at \url{https://ssd.jpl.nasa.gov/tools/periodic_orbits.html}, an accompanying explanation at \cite{Vaquero2018Poincare:Tool}, and the specific ones we use in our code repository.

\subsection{Discovery of local Poincar\'e Maps of the 3BP}
\label{sec: RES,Disc local Poin}
To achieve stabilisation of a UPO we require a sufficiently accurate Poincar\'e map. In this section we verify our mapping discovery using the desired Lyapunov orbit. We outline the mappings discovered at various surfaces of section and how our augmentation strategy enhances mapping accuracy and influences the variance of the mapping coefficients. 

We begin by discovering a Poincar\'e map for the desired UPO from the Lyapunov orbit family, following the approach outlined in \S \ref{sec: Learning maps}. First, we consider the initial conditions of the desired UPO $\mathbf{\Bar{x}}$ and $m=10$ neighbouring UPOs, assembling these in $\mathbf{\Bar{X}}$. As discussed in \S \ref{sec: Learning maps}, these neighbouring UPOs are selected from the orbit catalogue \cite{Vaquero2018Poincare:Tool} using regular intervals of the Jacobi constant $C$, specifically $1.75\times10^{-4}$ LU$^2$/TU$^2$ for Lyapunov orbits. 
For now, we omit the augmentation step in \S \ref{sec: Learning maps} and simulate these trajectories, recording two crossings of surface of section $\Sigma_{1,L}$ \eqref{eq:S1}. We utilise SINDy, as detailed in \S \ref{sec: Learning maps}, to discover a nonlinear Poincar\'e map and then linearise around $\mathbf{\Bar{x}}$ to obtain
\begin{equation}
\mathbf{A}_{\Sigma_{1,L} \text{ - no aug}} = \left[ \begin{array}{cccccc}
0 & 0 & 0 & 28.867 & -0.26859 & 0 \\
0 & 0 & 0 &  0 & 0 & 0 \\
0 & 0 & 0 & 0 & 0 & 0 \\
0 & 0 & 0 & 220.20 & 0 & 0 \\
0 & 0 & 0 & -97.934 & 1.0000 & 0 \\
0 & 0 & 0 & 0 & 0 & 0 \\
\end{array} \right],
\end{equation}
which satisfies \eqref{eq: lin sys}. If we use the augmentation strategy shown in Figure \ref{fig: aug diag} with a perturbation of $\delta v = 2.5 \times 10^{-7}$ LU/TU this leads to a very different map
\begin{equation}
\mathbf{A}_{\Sigma_1,L \text{ - aug}} = \left[ \begin{array}{cccccc}
1231.4 & 0 & 0 &   14.450 &  330.47 & 0 \\
       0  & 0 & 0 &    0    &    0    & 0 \\
       0  & 0 & 0 &    0    &    0    & 0 \\
  9431.8 & 0 & 0 &  109.78 & 2533.3 & 0 \\
 -4175.9 & 0 & 0 &  -49.044 & -1120.6 & 0 \\
       0  & 0 & 0 &    0    &    0    & 0
\end{array} \right],
\label{eq:Lyap A no aug}
\end{equation}
which since the first column is no longer zero, describes how a deviation from $\mathbf{\Bar{x}}$ grows due to changes in $x$.

To investigate the accuracy of these mappings, we compute the monodromy matrix $\mathbf{M}$ \eqref{eq: Monodromy} by integrating the system of 42 differential equations described in \S \ref{sec: Invariant M} over one period, starting from $\mathbf{\Bar{x}}$. 
We know the eigenvalues of $\mathbf{M}$ should match those of the linearised Poincar\'e map $\mathbf{A}$ \cite{Koon2011DynamicalDesign.}, and therefore we compare the eigenvalues of $\mathbf{A}$ in Table \ref{tab: eigs L S1}, with and without the use of the augmentation strategy and various perturbation magnitudes $\delta v$, to those of $\mathbf{M}$. 
As discussed in \S \ref{sec: Lyap fam}, the unstable and stable eigenvalues of $\mathbf{M}$ are the reciprocal pair $\lambda_u = 219.5457$ and $\lambda_s = 0.0046$, thus we compare the total absolute error, i.e. sum of the absolute differences between the respective eigenvalues $\sum |\lambda_{\mathbf{M},i} - \lambda_{\mathbf{A},i}|$, rather than relative error, since $\lambda_s << 1$ and we aim to weight the errors equally. Notably, for any model discovered with a zero column this will always give an eigenvalue of zero close to $\lambda_s = 0.0046$, but this provides no information on the stable manifold.
Table \ref{tab: eigs L S1} shows that the total error drops significantly when using our augmentation strategy, up to three orders of magnitude at the minimum when a perturbation of $\delta v = 2.5\times10^{-7}$ LU/TU. Importantly this is achieved in a sample-efficient manner utilising only 55 data points in $\mathbf{X}_1$ and corresponding ones in $\mathbf{X}_2$, where each data point is a different crossing of $\Sigma$. Additionally, since the 3BP is a volume-preserving system, any flow map should have a determinant equal to one. Table \ref{tab: eigs L S1} highlights that a smaller total eigenvalue error corresponds to the $\det(\mathbf{A)}$ approaching one, which is expected since the determinant equals the product of the eigenvalues. These properties indicate that the magnitude of $\delta v$ plays a critical role in improving the accuracy and near volume-preserving behaviour of the map.
\begin{table}
    \centering
    \caption{Comparison of eigenvalues for linearised Poincar\'e map $\mathbf{A}$ at $\Sigma_{1,L}$ to those of $\mathbf{M}$, for various magnitudes of augmentation perturbation $\delta v$, showing up to a three orders of magnitude reduction in total eigenvalue error with our augmentation strategy. Final column shows $\det(\mathbf{A})$ approaches unity for the same $\delta v$, indicating the map is approaching a volume-preserving structure.}
    % \begin{tabular}{l cccc}
    %     \toprule
    %     & \multicolumn{4}{c}{Eigenvalues} \\
    %     \cmidrule(lr){2-5}
    %     \textbf{} & \textbf{Unstable} & \textbf{Neutral} & \textbf{Stable} & \textbf{Total Error}\\
    %     \midrule
    %     $\mathbf{M}$ & 219.5457 & 1.0000 & 0.0046 & $\cdot$\\
    %     $\mathbf{A}_{\Sigma_{1,L}} -\text{no augmentation}$ & 220.3278 & 1.0000 & 0.0000 & 0.7868\\
    %     $\mathbf{A}_{\Sigma_{1,L}} - \delta v = 2.5\times10^{-9} $ LU/TU & 219.5400 & 1.0020 & 0.0006 & 0.0117\\
    %     $\mathbf{A}_{\Sigma_{1,L}} - \delta v = 2.5\times10^{-7} $ LU/TU & 219.5454 & 1.0004 & 0.0044 & 0.0008\\
    %     $\mathbf{A}_{\Sigma_{1,L}} - \delta v = 2.5\times10^{-6} $ LU/TU & 219.5408 & 0.9998 & 0.0030 & 0.0067\\
    %     \bottomrule
    % \end{tabular}
    \begin{tabular}{l ccccc}
    \toprule
    & \multicolumn{5}{c}{Eigenvalues \& det(A)} \\
    \cmidrule(lr){2-6}
    & \textbf{Unstable}
    & \textbf{Neutral}
    & \textbf{Stable}
    & \textbf{Total Error}
    & \text{det(\textbf{A})}\\
    \midrule
    $\mathbf{M}$
      & 219.5457
      & 1.0000
      & 0.0046
      & $\cdot$
      & 1.0000\\
    $\mathbf{A}_{\Sigma_{1,L}} -\text{no augmentation}$
      & 220.3278
      & 1.0000
      & 0.0000
      & 0.7868
      & 0.0000 \\
    $\mathbf{A}_{\Sigma_{1,L}} - \delta v = 2.5\times10^{-9}$ LU/TU
      & 219.5400
      & 1.0020
      & 0.0006
      & 0.0117
      & 0.1386 \\
    $\mathbf{A}_{\Sigma_{1,L}} - \delta v = 2.5\times10^{-7}$ LU/TU
      & 219.5454
      & 1.0004
      & 0.0044
      & 0.0008
      & 0.9741 \\
    $\mathbf{A}_{\Sigma_{1,L}} - \delta v = 2.5\times10^{-6}$ LU/TU
      & 219.5408
      & 0.9998
      & 0.0030
      & 0.0067
      & 0.6688 \\
    \bottomrule
\end{tabular}
    \label{tab: eigs L S1}
\end{table}

\begin{figure}
    \centering
    \includegraphics[width=0.5\linewidth]{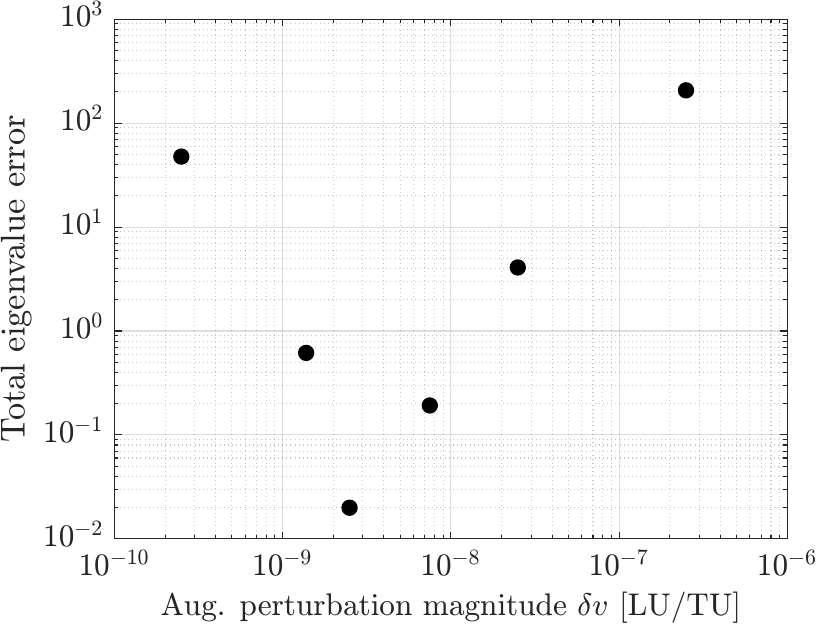}
    \caption{Total eigenvalue error of $\mathbf{A}$ at $\Sigma_{2,L}$ versus augmentation perturbation magnitude $\delta v$, highlighting a decrease to a minimum error of 0.02 at $\delta v = 2.50\times10^{-9}$ LU/TU before rising again}
    \label{fig: S2 err}
\end{figure}

\begin{figure}
  \centering
  \begin{subfigure}[b]{0.48\textwidth}
    \includegraphics[width=\linewidth]{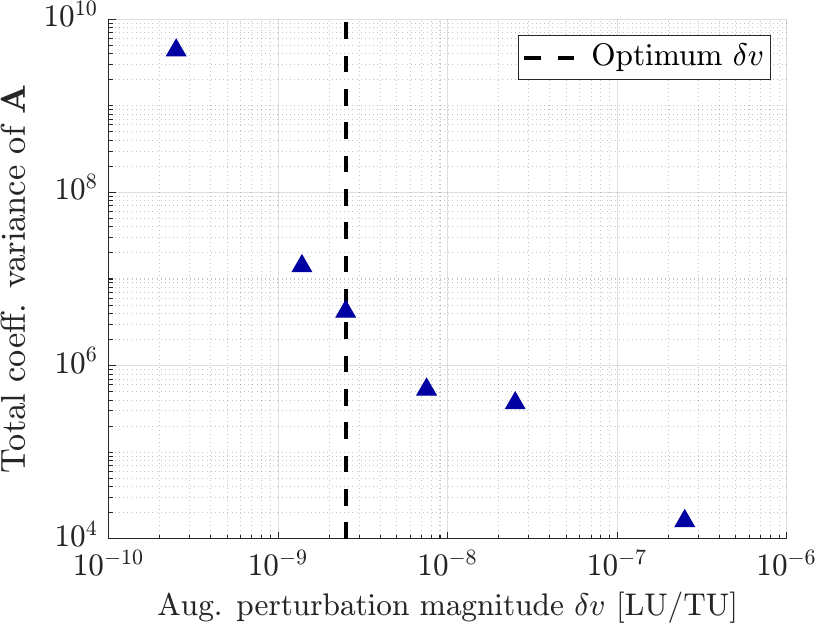}
    \caption{A near logarithmic decrease in total coefficient variance of $\mathbf{A}$ as augmentation perturbation magnitude $\delta v$ increases.}
    \label{fig: SL2 coef var}
  \end{subfigure}
  \hfill
  \begin{subfigure}[b]{0.48\textwidth}
    \includegraphics[width=\linewidth]{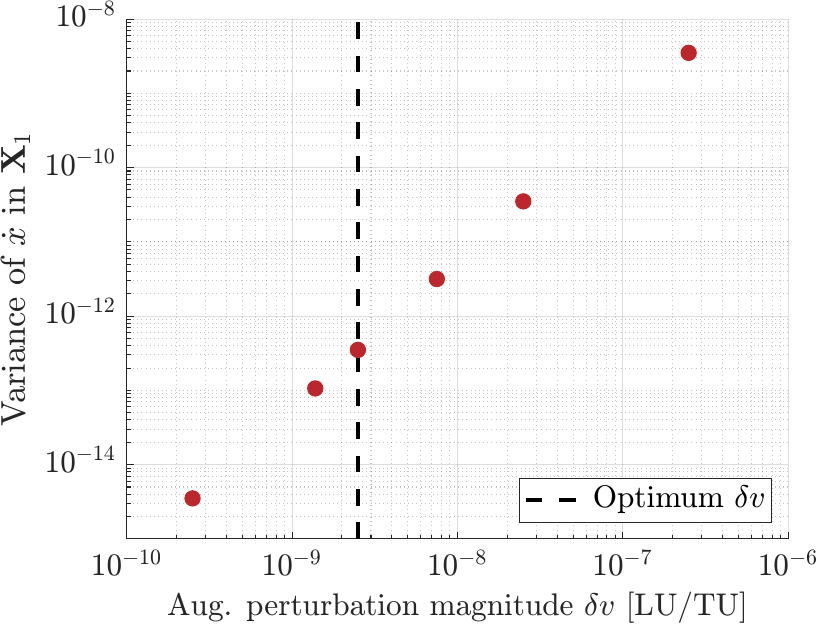}
    \caption{Variance of $\Dot{x}$ in training data $\mathbf{X}_1$ increasing almost logarithmically as augmentation perturbation magnitude $\delta v$ increases.}
    \label{fig: SL2 xdot var}
  \end{subfigure}
  \caption{Inverse trends in total coefficient variance of $\mathbf{A}$ and variance in the training data as $\delta v$ increases when discovering a Poincar\'e map at $\Sigma_{2,L}$}
  \label{fig:SL2}
\end{figure}

We now move our attention to discovering a mapping at the alternate surface of section $\Sigma_{2,L}$ \eqref{eq:S2} and investigating the influence of $\delta v$ further. The eigenvalues of the linearised maps $\mathbf{A}$ and $\mathbf{M}$ remain the same at both surfaces of section since they describe the dynamics of the same UPO and the maps at either surface are conjugate.
Again, we vary the augmentation perturbation $\delta v$ and record the total eigenvalue error of the linearised discovered mapping $\mathbf{A}$, independent of the number of data points, which remains 55 for all $\delta v$ investigated at $\Sigma_{2,L}$. Figure \ref{fig: S2 err} shows how the total eigenvalue error of $\mathbf{A}$ varies with $\delta v$, initially decreasing by three orders magnitude to a minimum of 0.02 at $\delta v = 2.50\times10^{-9}$ LU/TU before sharply rising again, indicating again that $\delta v$ must be adjusted to discover the most accurate map.
% To investigate the effect of the augmentation strategy on the accuracy of $\mathbf{A}$, We vary the augmentation perturbation magnitude $\delta v$, and plot this versus the total eigenvalue error of the discovered mapping in Figure \ref{fig: S2 err}. 
% Increasing the perturbation from $\delta v = 2.50\times10^{-10}$, the total eigenvalue error decreases by three orders magnitude to a minimum of 0.02 at $\delta v = 2.50\times10^{-9}$ before sharply rising again. 
This minimum corresponds to the linearised map
\begin{equation}
    \mathbf{A}_{\Sigma_{2,L}} = \left[ \begin{array}{cccccc}
1204.3 & 0 & 0 & 3.0409 & 9.0089 & 0 \\
0.0027936 & 0 & 0 & 0 & 0 & 0 \\
0 & 0 & 0 & 0 & 0 & 0 \\
43827 & 0 & 0 & 109.76 & 328.13 & 0 \\
-146190 & 0 & 0 & -369.43 & -1093.5 & 0 \\
0 & 0 & 0 & 0 & 0 & 0 \\
\end{array} \right].
\end{equation}
We investigate this behaviour further by using the ensemble methods described in \S \ref{sec: Learning maps}, to measure the variance of the model coefficients. We choose a sparsity promoting parameter of $\lambda_{\text{sparse}}=$0.02, $N_e$ = 100 ensemble models, and set the threshold inclusion probability to $\rho=0.6$. These parameters are effective for our study and we direct the reader to the accompanying code to explore this further.
Denoting each coefficient in $\mathbf{A}$ as $\xi_{A,i}$, we calculate the total variance of these coefficients $\sum \text{Var}(\xi_{A,i})$ across the 100 models. Figure \ref{fig: SL2 coef var} shows that total variance of these model coefficients decreases almost logarithmically as $\delta v$ increases. 
The variance in the training data $\mathbf{X}_1$ remains the same order of magnitude for all the states apart from $\Dot{x}$, which shows an almost logarithmic increase with $\delta v$ shown in Figure \ref{fig: SL2 xdot var}.
The magnitude of $\delta v$ has a significant influence on the accuracy of the mapping discovery, demonstrated through these trends at both surfaces of section, and that there is a compromise to be made between the variance in the training data and in the model coefficients.
The same analysis as Figures \ref{fig: S2 err}, \ref{fig: SL2 coef var}, \ref{fig: SL2 xdot var} can be performed for $\Sigma_{1,L}$ but this has been omitted for brevity, as similar trends are observed. 

Importantly, the most accurate $\mathbf{A}$ at $\Sigma_{1,L}$ has a total eigenvalue error 60\% less than for the most accurate at $\Sigma_{2,L}$, suggesting it is easier to discover a more accurate model at $\Sigma_{1,L}$. 
This is likely due to the Lyapunov orbit family passing very close to the secondary mass, exhibiting a significantly increased sensitivity to perturbations  in the region examined in $\Sigma_2$. This is reflected in significant curvature of the orbits close to $x = 1 - \mu$ in Figure \ref{fig:Lyap Fam}.
Recalling our measure of sensitivity in \S \ref{sec: Invariant M} we take the norm of the monodromy matrix and find $||\mathbf{M}|| = 5.9\times10^4$ at $\Sigma_{1,L}$, an order of magnitude smaller than $||\mathbf{M}|| = 5.7\times10^5$ at $\Sigma_{2,L}$. This suggests that a more accurate map may be discovered at the surface of section where the sensitivity to small perturbations is less, i.e where $||\mathbf{M}||$ is smaller.

\subsection{Stabilisation of Lyapunov UPOs}
\begin{figure}
    \centering
    \includegraphics[width=1.0\linewidth]{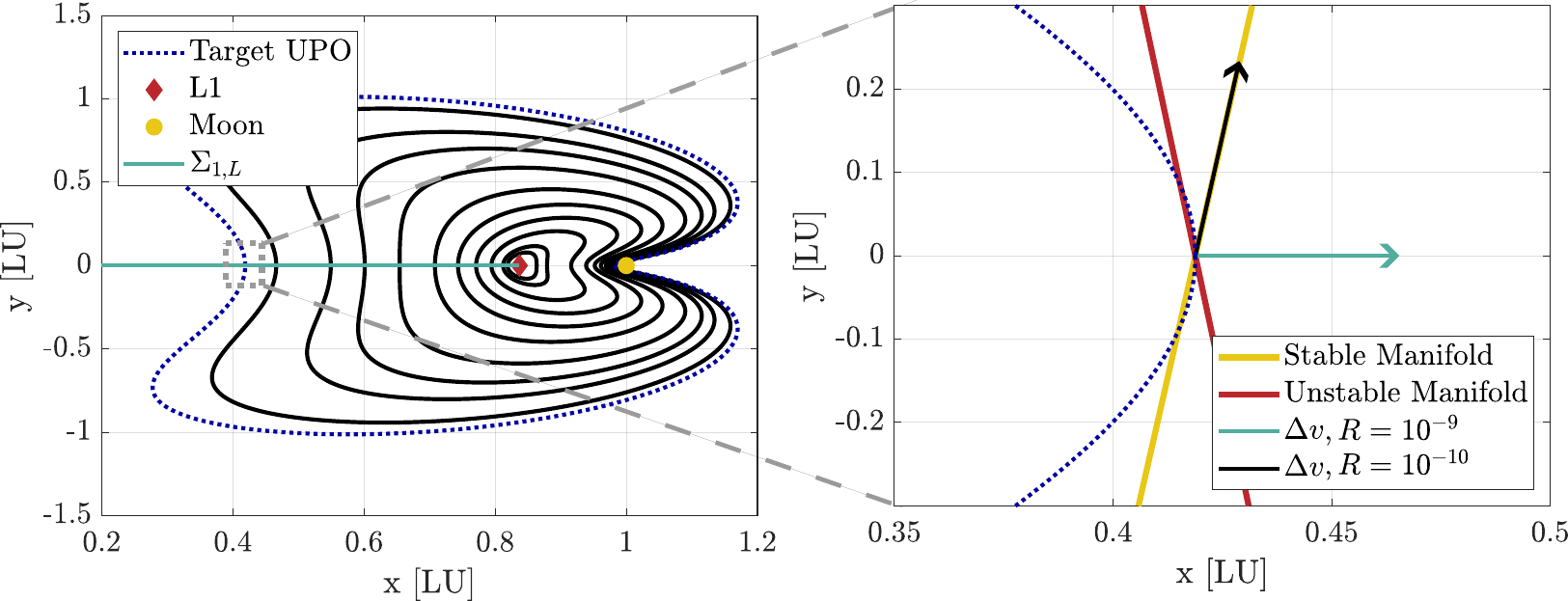}
    \caption{The stabilised UPO within the Lyapunov family (left), showing its intersection with the surface of section $\Sigma_{1,L}$. The right panel zooms in at $\Sigma_1$ illustrating how the reduction in feasibility radius $R$ reorients the control impulse $\Delta v$ to act along the local stable manifold.}
    \label{fig:Manifolds}
\end{figure}

\begin{figure}
  \centering
  \begin{subfigure}[b]{\textwidth}
    \begin{minipage}{0.7\linewidth}
        \centering
        \includegraphics[width=\linewidth]{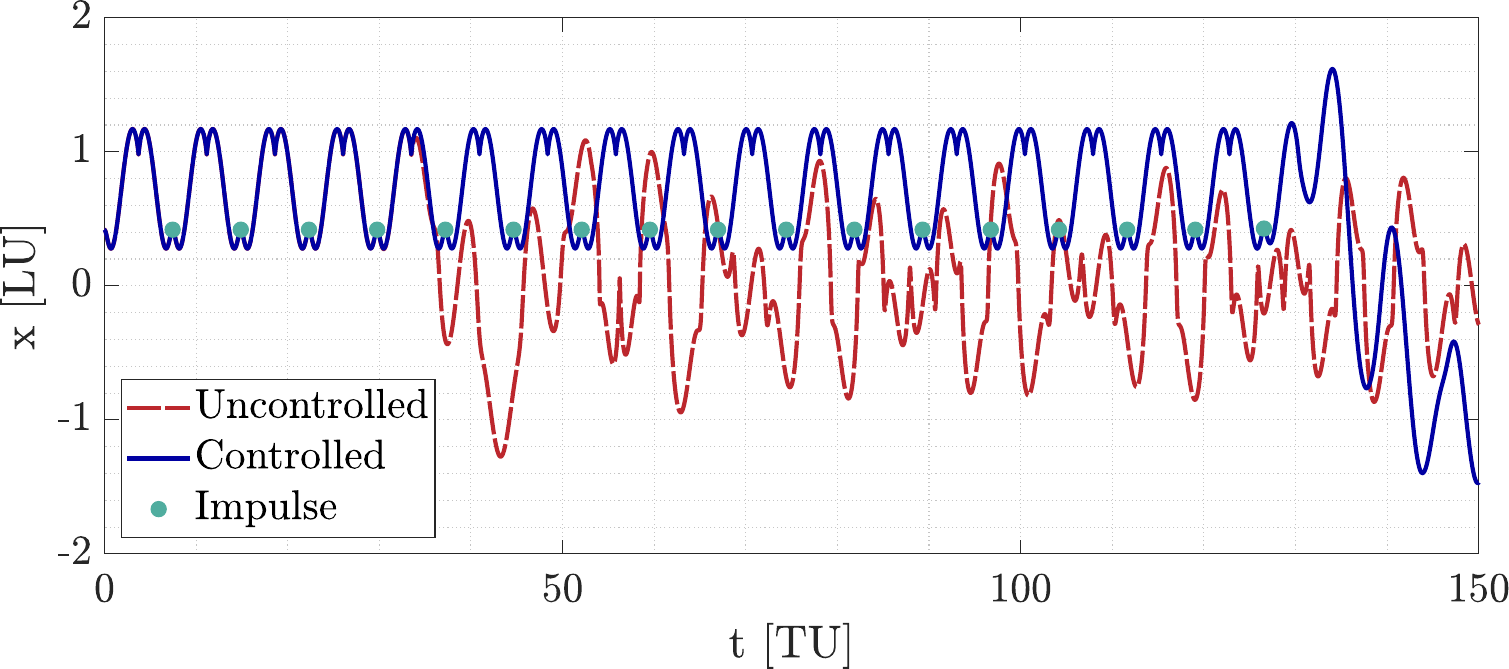}
    \end{minipage}
    \hfill
    \begin{minipage}{0.25\linewidth}
        \centering
        \includegraphics[width=\linewidth]{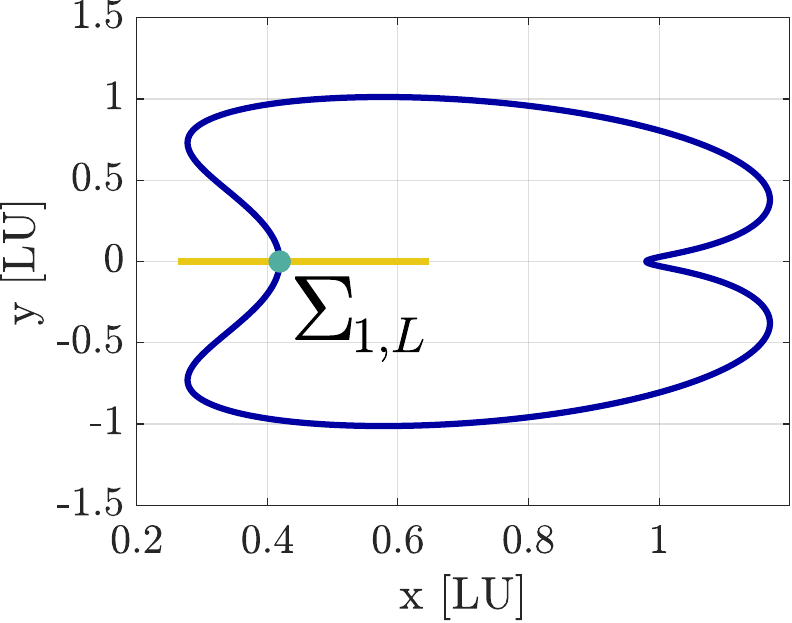}
        \begin{equation*}
            R =10^{-10} \quad\theta = 0.01 \degree
        \end{equation*}
        \begin{equation*}
        \Delta v_{1:14} = 9.25 \times 10^{-8} \text{ m/s}
        \end{equation*}
    \end{minipage}
    \caption{Behaviour and properties for best controller at $\Sigma_{1,L}$}
    \label{fig:x vs t aug sec1}
  \end{subfigure}
  
  \begin{subfigure}[b]{\textwidth}
    \begin{minipage}{0.7\linewidth}
        \centering
        \includegraphics[width=\linewidth]{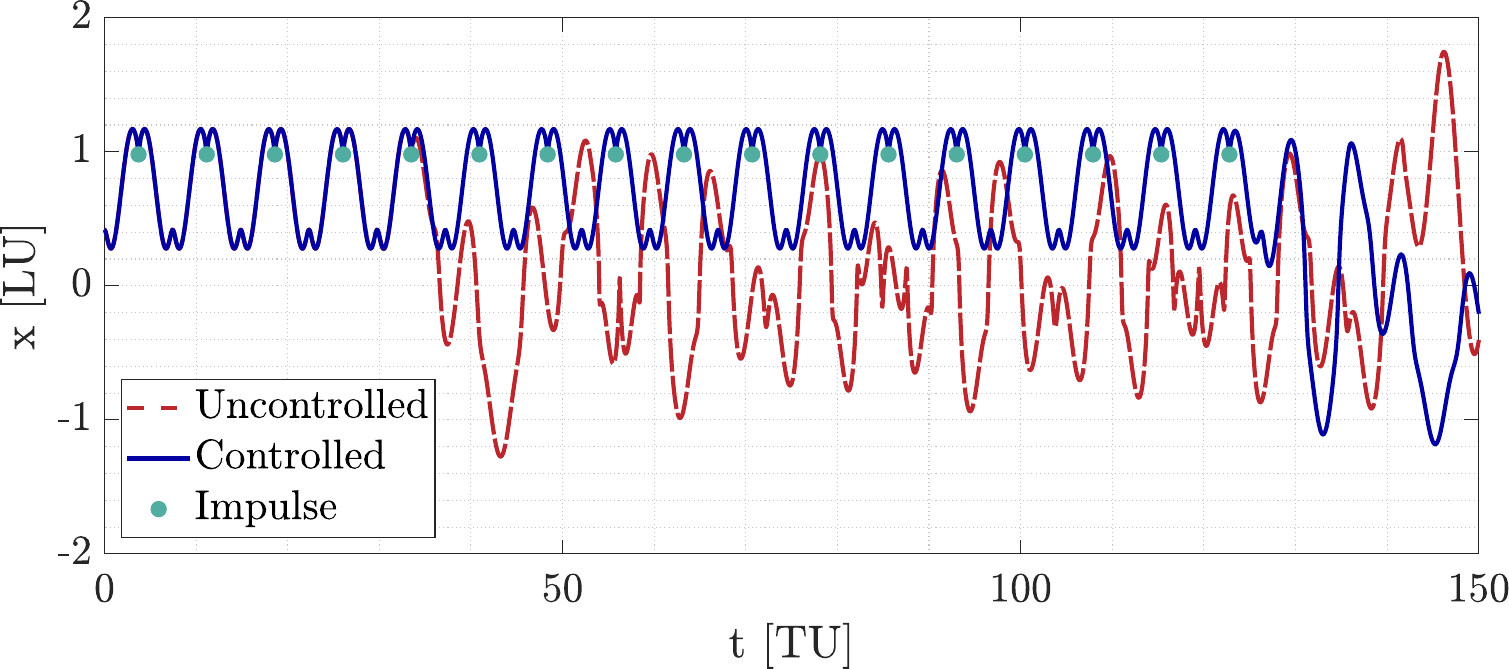}
    \end{minipage}
    \hfill
    \begin{minipage}{0.25\linewidth}
        \centering
        \includegraphics[width=\linewidth]{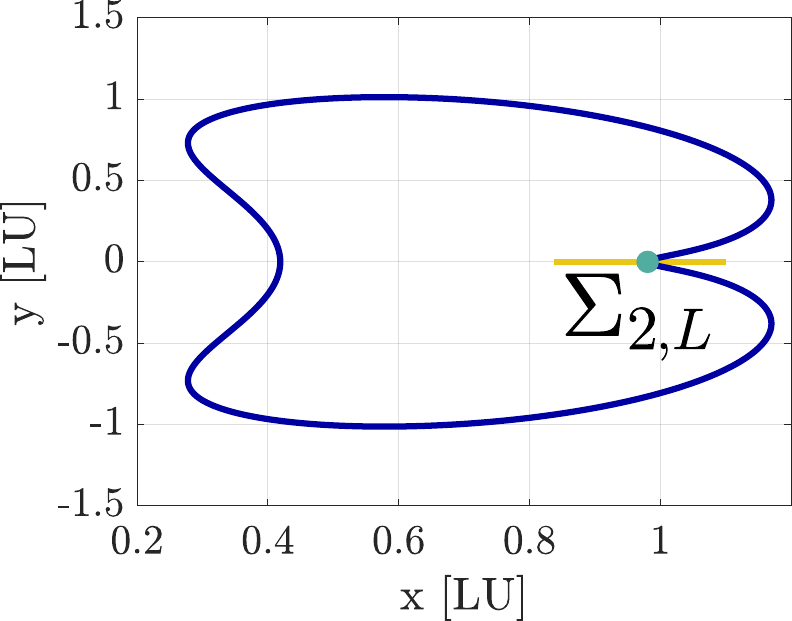}
        \begin{equation*}
            R = 1.19 \times 10^{-10} \quad  \theta = 0.27 \degree
        \end{equation*}
        \begin{equation*}
        \Delta v_{1:14} =  3.33 \times 10^{-5} \text{ m/s}
        \end{equation*}
    \end{minipage}
  \caption{Behaviour and properties for best controller at $\Sigma_{2,L}$}
  \label{fig:x vs t aug sec2}
  \end{subfigure}
  \caption{Left) Time evolution of the $x$ component showing the Lyapunov UPO stabilised with periodic impulses acting at each $\Sigma$. Right) The intersection of the UPO with each $\Sigma$ alongside the key controller properties, showing a smaller angle between the local stable manifold and the impulse and reduced control cost at $\Sigma_{1,L}$.}
\end{figure}

Once a sufficiently accurate mapping has been found we may solve the LMI problem \eqref{eq: LMI form} subject to the constraint \eqref{eq: R}, where 
\begin{equation}
    \mathbf{R} = 
\begin{bmatrix}
0 & 0 & 0 & 1 & 0 & 0 \\
0 & 0 & 0 & 0 & 1 & 0 \\
0 & 0 & 0 & 0 & 0 & 1
\end{bmatrix}^T,
\end{equation}
allowing only the velocity states may be altered. First, we solve the problem using $\mathbf{A}_{\Sigma_1,L}$ discovered with $\delta v = 2.5\times 10^{-7}$ LU/TU at $\Sigma_{1,L}$, and setting the feasibility radius to $R = 10^9$ which in effect is not active as a constraint, results in the control matrix
\begin{equation}
\mathbf{K} = \left[ \begin{array}{cccccc}
-7.6699 & 0 & 0 & -1.0000 & 0 & 0 \\
-3.3908 & 0 & 0 & 0 & -1.0000 & 0 \\
0 & 0 & 0 & 0 & 0 & 0 \\
\end{array} \right].
\label{res: K1a}
\end{equation}
Using this controller led to the successful stabilisation of the selected Lyapunov orbit for 17 periods at $\Sigma_{1,L}$, which can be seen on the left panel of Figure \ref{fig:x vs t aug sec1}. 
We decrease the value of $R$ and compare the effectiveness of each controller by calculating the total impulse over the first 14 periods $\Delta v_{1:14}$ as a representative measure. The units are re-dimensionalised so the reader may have some intuition about how small these impulses really are. Table \ref{tab: Controller Prop S1} shows a gradually decreasing magnitude in $\Delta v_{1:14}$ initially, and then a sudden reduction by two orders of magnitude between $R = 10^{-9}$ and $10^{-10}$. The controller corresponding to $R=10^{-10}$ is of the form
\begin{equation}
    \mathbf{K} = \left[ \begin{array}{cccccc}
-0.1613 & 0 & 0 & -0.0019 & -0.0433 & 0 \\
-3.7156 & 0 & 0 & -0.0433 & -0.9978 & 0 \\
0 & 0 & 0 & 0 & 0 & 0 \\
\end{array} \right].
\end{equation}
Decreasing $R$ beyond this did not result in improved properties.
We investigate this sudden reduction in control cost by computing the local stable and unstable manifold directions from the monodromy matrix $\mathbf{M}$ we calculated in Table \ref{tab: eigs L S1}, and comparing them to impulse directions. For the unstable eigenvalue $\lambda_u = 219.5457$ the associated eigenvector is $\nu_{u} = [-0.0217,0.5005, 0.0000, -0.8623, 0.0736, 0.0000]$ and for the stable eigenvalue $\lambda_s = 0.0046$ the associated eigenvector is $\nu_s = [-0.0217,-0.5005, 0.0000, 0.8623, 0.0736, 0.0000]$. We extract the $x$ and $y$ components into one vector and normalize it to obtain a direction vector. Similarly, we normalize the velocity impulse directions, which all follow the same direction each crossing of $\Sigma_2$, and compute the angle between the impulse and local stable manifold $\nu_s$ denoted  $\theta$. From Table \ref{tab: Controller Prop S1}, we observe that $\theta$ decreases $5 \degree$  with the initial reduction in $R$ but then suddenly drops $82.57 \degree$ to $0.01 \degree$ between $R = 10^{-9}$ and $10^{-10}$, visualised in Figure \ref{fig:Manifolds}. This corresponds with the large reduction in total impulse, indicating that impulses directed along the local stable manifold minimise control activity. 
\begin{table}
    \centering
    \caption{The effect of reducing the feasibility radius $R$ on the controller behaviour at $\Sigma_{1,L}$: the eigenvalues magnitude gradually increases and 
    the re-dimensionalised total impulse over the first 14 periods $\Delta v_{1:14}$ suddenly decreases between $R=10^{-9}$ and $R=10^{-10}$, corresponding to the impulse reorienting to align with the local stable manifold $\nu_s$, with the angle between them $\theta$ dropping to near zero.}
    \begin{tabular}{cccc}  
        \toprule
        $R$ &  \shortstack{$\theta [\degree]$} & $\Delta v_{1:14}$ [m/s] &  eig($\mathbf{A+ARK}$) \\
        \midrule
        $10^9$ & 87.52 & 7.23$\times10^{-6}$ &  7.39$\times10^{-4}$, 8.37$\times10^{-10}$,  -6.62$\times10^{-14}$ \\
        \midrule
        $10^{-7}$ & 82.39 & 3.30$\times10^{-6}$ & 
        0.393, \shortstack{-1.69$\times10^{-5}$ + 7.95$\times10^{-4}$i, -1.69$\times10^{-5}$ - 7.95$\times10^{-4}$i } \\
        \midrule
        $10^{-9}$ & 82.62 & 3.52$\times10^{-6}$ & 0.390, 3.71$\times 10^{-5}$+ 1.44$\times 10^{-3}$i, 3.71$\times 10^{-5}$ - 1.44$\times 10^{-3}$i \\
        \midrule
        $10^{-10}$ & 0.01 & 9.25$\times10^{-8}$ & 0.279 + 0.356i, 0.279 - 0.356i, 1.605$\times 10^{-3}$\\
        \bottomrule
    \end{tabular}
    \label{tab: Controller Prop S1}
\end{table}

Using the same procedure, we analyse the controller behaviour at $\Sigma_{2,L}$. In Table \ref{tab: Cont prop S2} we again observe that the eigenvalues increase in magnitude and there is a gradual reduction in the control cost, but no sudden decrease as seen in Table \ref{tab: Controller Prop S1}. Despite this, at $R = 1.19\times10^{-10}$, the smallest possible value that ensures the eigenvalues $\lambda$ are such that $|\lambda| \leq 1$, the associated controller
\begin{equation}
    \mathbf{K} = \begin{bmatrix}
    -40.5032      & 0  & 0 & -0.1029 & -0.3028 & 0 \\
    -119.9963     & 0 & 0 & -0.3028 & -0.8977 & 0 \\
    0             & 0                     & 0 & 0       & 0       & 0
    \end{bmatrix},
\end{equation}
further reduces control cost and the impulse is again reoriented to align with the local stable manifold at an angle of just $0.27 \degree$. Importantly the impulses along the local stable manifold are also used to minimise cost in station-keeping strategies used in the NASA ARTEMIS mission \cite{Folta2010StationkeepingARTEMIS} and other methods proposed after \cite{Pavlak2012StrategySystem}. Constraining the magnitude of the decision variables $\mathbf{Q}$ and $\mathbf{Y}$ in \eqref{eq: LMI form} causes an increase in the magnitude of the eigenvalues shown in Tables \ref{tab: Controller Prop S1} and \ref{tab: Cont prop S2}, indicating a less aggressive controller. Critically, this forces impulses to align with the local stable manifold, the locally-optimal path since this direction is associated with the decay of perturbations and a stable eigenvalue.

Figures \ref{fig:x vs t aug sec1} and \ref{fig:x vs t aug sec2} show the time evolution of the stabilised UPOs for the best controllers at $\Sigma_1$ and $\Sigma_2$ respectively, alongside the corresponding surface of section. The feasibility radius $R$ is very similar but the angle between the impulse and local stable manifold $\theta$ is an order of magnitude less at $\Sigma_1$, and the total re-dimensionalised impulse $\Delta v_{1:14}$ is three orders of magnitude less.

\begin{table}
    \centering
    \caption{
    The effect of reducing the feasibility radius $R$ on the controller behaviour at $\Sigma_{2,L}$: the eigenvalues magnitude gradually increases and the re-dimensionalised total impulse over the first 14 periods $\Delta v_{1:14}$ gradually decreases. At the smallest $R$ the impulse $\delta v$ again acts along the local stable manifold $\nu_s$, with the angle $\theta$ between them almost zero.
    }
    \begin{tabular}{cccc}  
        \toprule
         $R$ & \shortstack{$\theta [\degree]$} & $\Delta v_{1:14}$ [m/s] & eig($\mathbf{A+ARK}$) \\
        \midrule
        $10^9$ & 71.58 & 2.24$\times10^{-4}$ & -0.00280, -7.07$\times10^{-14}$,  5.30$\times10^{-12}$ \\
        \midrule
        $10^{-7}$ & 51.88 & 8.2$\times10^{-5}$ & 
        -0.456, 0.00700, -0.00281\\
        \midrule
        $10^{-9}$ & 48.63 & 7.50$\times10^{-5}$ & -0.504, -0.00160, -0.00131 \\
        \midrule
        1.19 $\times 10^{-10}$ & 0.27 & 3.33$\times 10^{-5}$ & -0.775, -0.0117 + 0.0137i, -0.0117 - 0.0137i\\
        \bottomrule
    \end{tabular}
    \label{tab: Cont prop S2}
\end{table}

\subsection{Stabilisation of Halo UPOs}
In this subsection we successfully demonstrate the stabilisation of a UPO from the halo orbit family, described in \S \ref{sec: Halo fam}, whose motion is not limited to the orbital plane. We utilise the same augmentation strategy, and found that changing the magnitude of the augmentation perturbation $\delta v$ does not result in significantly different mappings being discovered. We omit further discussion of this for brevity since our focus is to demonstrate our strategy can stabilise non-planar UPOs. First, we consider the surface of section $\Sigma_{1,H}$ \eqref{eq:H1} with a perturbation of $\delta v = 2.5 \times10^{-9}$ LU/TU, which leads to the discovered and subsequently linearised map 
\begin{equation}
    \mathbf{A}_{\Sigma_{1,H}} = 
\begin{bmatrix}
  59.621   & 0 & 219.65  & 169.59  &  -38.570  & 157.14 \\
   0      & 0 &   0     &   0     &    0     &   0    \\
 -125.73  & 0 & 116.20  & 157.14  & -177.12  & 109.20 \\
 -868.37  & 0 & -357.15 &  60.05  & -712.36  & -123.57\\
  191.49  & 0 & -178.76 & -242.25 &  272.21  & -168.74\\
  737.25  & 0 &  675.12 &  278.59 &  439.97  &  386.98
\end{bmatrix}.
\end{equation}
The eigenvalues of this map \{447.0 + 175.63i, 447.0 - 175.63i, -0.0027 + 0.0047i, -0.0027 - 0.0047i, 1.00\}, do not follow the form \eqref{eq: halo eigs} and similarly the eigenvalues of the computed monodromy matrix contain two complex conjugate pairs that are reciprocals. The magnitude of the pairs matches the expected values of the eigenvalues $\lambda_u = 480.2979$ and $\lambda_s = 0.0021$ however, and the imaginary components likely appear due to a loss of symplectic structure in the integration scheme, even though we use a high order integrator at a low tolerance. Nevertheless, the magnitudes of the eigenvalues remain accurate and this proves sufficient for control purposes.

\begin{figure}
  \centering
  \begin{subfigure}[b]{0.6\textwidth}
    \includegraphics[width=\linewidth]{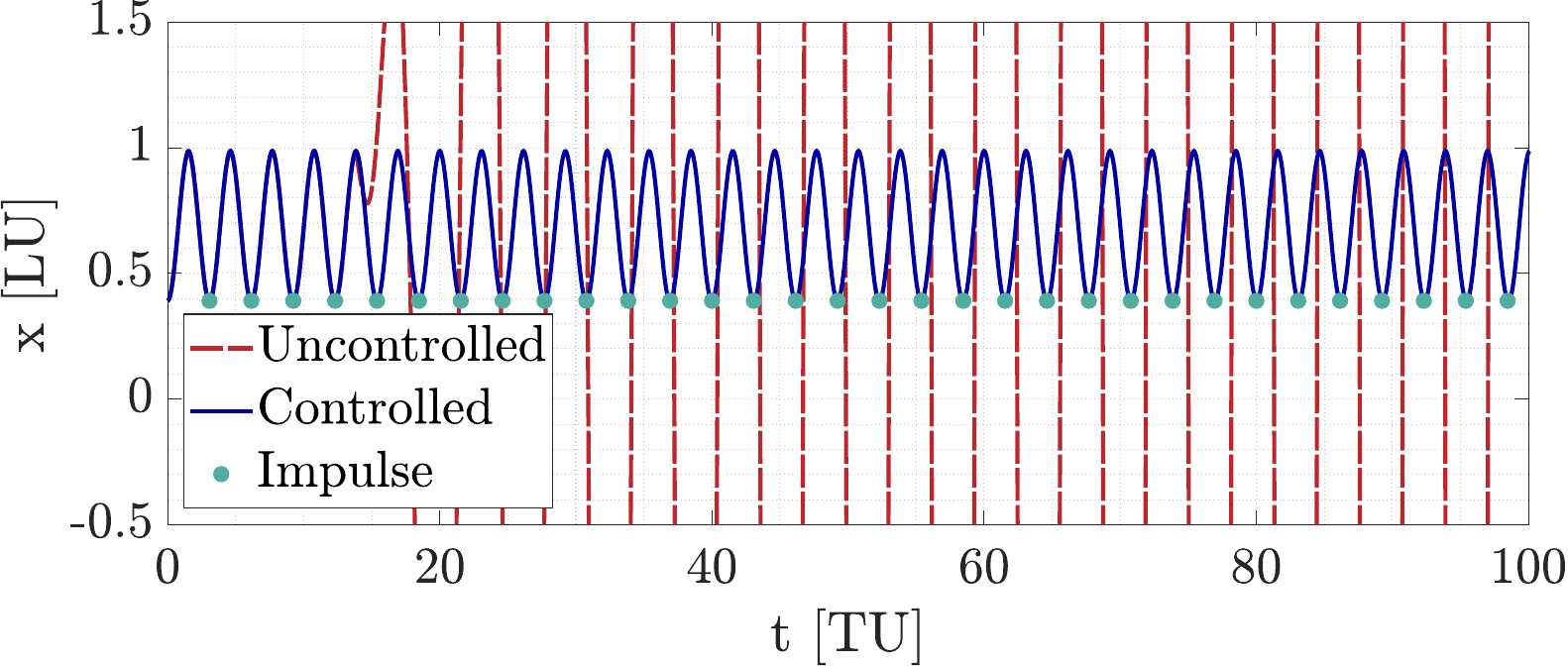}
    \caption{The $x$ component of the stabilised halo UPO for 30 periods, showing the recurrent periodic behaviour maintained by the controller applying an impulse once every period.}
    \label{fig: x vs t halo}
  \end{subfigure}
  \hfill
  \begin{subfigure}[b]{0.35\textwidth}
    \includegraphics[width=\linewidth]{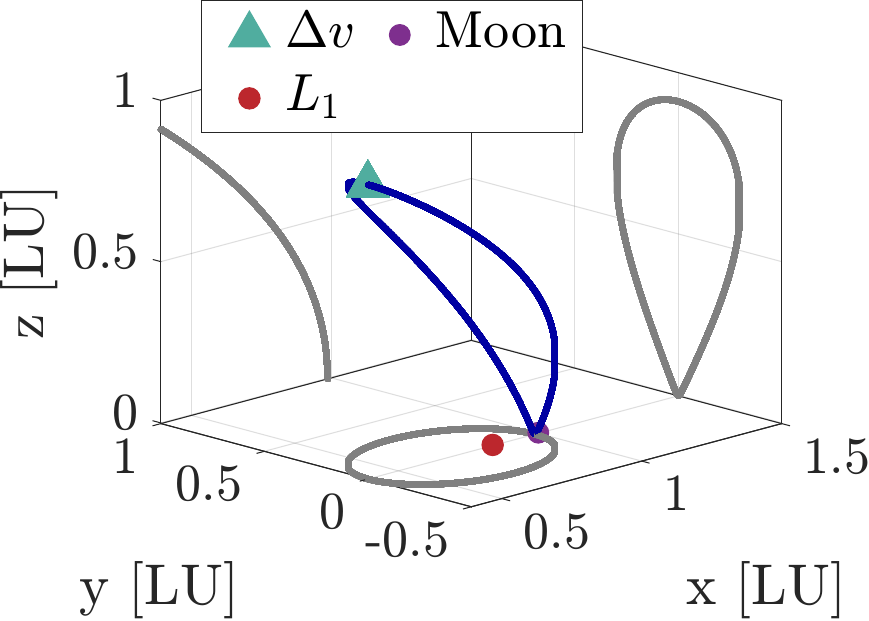}
    \caption{Stabilised halo UPO with projections onto each plane with the impulse point at the top of the orbit shown at $\Sigma_{1,H}$.}
    \label{fig: 3D halo}
  \end{subfigure}
  \caption{Visualisation of the stabilised UPO from the halo orbit family.}
  \label{fig: Halo stabilised}
\end{figure}

We use the map $\mathbf{A}_{\Sigma_{1,H}}$ to stabilise the UPO, through solving the same LMI problem as before \eqref{eq: LMI form}. Reducing $R$ only results in a small reduction in control cost. Taking the minimum feasibility radius $R=9.5\times10^{-10}$ for which the eigenvalues $\lambda$ of the controlled system $\mathbf{A+BK}$ are such that $|\lambda| \leq 1$, this results in the controller
\begin{equation}
    \mathbf{K} = 
\begin{bmatrix}
 -0.1762 & 0      & -0.6771 & -0.5243 &  0.1244 & -0.4842 \\
 -1.1521 & 0      & -0.4238 &  0.1244 & -0.9677 & -0.1257 \\
 -0.4688 & 0      & -0.7751 & -0.4843 & -0.1257 & -0.5089
\end{bmatrix},
\end{equation}
which successfully stabilises the UPO for more than 90 periods. In Figure \ref{fig: x vs t halo} we show the $x$ component of the stabilised UPO, with only the first 32 periods (100 TU) for brevity, demonstrating the recurrent periodic behaviour maintained by the impulse once every period. The total impulse of the first 32 periods is $\Delta v_{1:33} = 3.09\times10^{-6}$ m/s. We illustrate the stabilised UPO in Figure \ref{fig: 3D halo}, highlighting the location of the impulse within $\Sigma_{1,H}$ at the maximum $z$ value.  

No meaningful mapping was discovered at $\Sigma_{2,H}$ \eqref{eq:H2} and consequently the UPO was not stabilised at this surface of section. Taking our measure of local sensitivity, the norm of the mondromy matrix at $\Sigma_{1,H}$ is $||\mathbf{M}|| = 1.4\times10^{3}$ and at  $\Sigma_{2,H}$ is $||\mathbf{M}|| = 2.3\times10^{6}$. This concurs with our results and suggestion in \S \ref{sec: RES,Disc local Poin} that a mapping is more easily discovered at a surface of section with a smaller magnitude of monodromy matrix. Intuitively, in Figure \ref{fig:Halo Fam} we see that at $\Sigma_{2,H}$ the orbits suddenly change direction around the moon, illustrating where this sensitivity arises from.

\section{Discussion}
\label{sec: Discussion}
In this paper, we have presented a sample-efficient and interpretable data-driven approach to achieve low-energy stabilisation of UPOs of the 3BP. We outline a novel augmentation strategy that utilises prior knowledge of UPOs and small velocity perturbations to collect section data, which enables the discovery of an accurate state-dependent local Poincar\'e map using the SINDy algorithm \cite{Brunton2015DiscoveringSystems} in a volme-preserving system. Our mapping discovery is effective in the low-data limit, utilising as few as 55 data points, in contrast to deep learning methods \cite{Bramburger2021DeepMappings} that  often require large volumes of data which is particularly costly to collect for chaotic systems. Importantly, our augmentation strategy enables us to adapt the OGY method \cite{Ott1990ControllingChaos} for state-dependent, rather than parameter-dependent control, exploiting the inherent sensitivity of chaos to achieve low-energy control. We utilise a series of LMIs to form a convex pole-placement problem which we solve using an interior point method \cite{Nesterov1994Interior-PointProgramming} in less than a second on a laptop. We add a constraint on the magnitude of the decision variables which reduces the control cost and results in the impulses acting along the local stable manifold.
Our work suggests that surfaces of section are best defined where the local dynamics are less sensitive to perturbations, which we quantify by the magnitude of the monodromy matrix and may serves as a guide to place the best surface.

First, we investigated the local mapping discovery on a selected UPO from the planar Lyapunov orbit family. Developing our novel augmentation strategy was motivated by the volume-preserving nature of conservative systems, meaning we cannot use a simple parameter sweep or random initialisations to sample the dynamics since the UPOs do not all fall on a strange attractor. The interpretable nature of SINDy enables us to verify the accuracy of the discovered mappings by comparing the eigenvalues of the linearised map $\mathbf{A}$ to those of the monodromy matrix $\mathbf{M}$ of the UPO. Subsequently, we found that the magnitude of the augmentation perturbation $\delta v$ has a significant effect on the accuracy and the variance of the model coefficients.
The most accurate map, with a total eigenvalue error two orders of magnitude less than at the alternate section, was discovered at the surface of section $\Sigma$ where the magnitude of the monodromy matrix $||\mathbf{M}||$ was the smallest, and by implication where the sensitivity to small perturbations is the least.

We successfully stabilise a UPO from the planar Lyapunov family at two different surfaces of section, and we found that constraining the magnitude of the decision variables in the LMI control formulation reduced the control cost. This constraint increases the magnitude of the controlled system's eigenvalues, indicating a less aggressive controller. This causes the impulses to act along the local stable manifold, the direction associated with a stable eigenvalue and decay in perturbation magnitude. Critically, this is similar to the original OGY method \cite{Ott1990ControllingChaos}, and this locally-optimal behaviour corresponds to existing control strategies for station-keeping in the 3BP \cite{Pavlak2012StrategySystem, Folta2009StationkeepingMission}. A key advantage of our approach is that we achieve this locally-optimal control behaviour in a computationally efficient manner, without the need for expensive direct optimisation methods.
We demonstrate that our method is easily applicable to non-planar orbits as well, successfully stabilising a UPO from the halo orbit family. In this case we only discover a meaningful map and achieve subsequent control of the UPO at the surface of section $\Sigma$ with the smaller magnitude of monodromy matrix. Furthermore, in the case of the Lyapunov UPO we found the control cost was three orders of magnitude less at the surface of section with a smaller value of $||\mathbf{M}||$. This further supports our suggestion that a surface of section is best placed, for mapping discovery and control, where $||\mathbf{M}||$ and thus the local sensitivity to small perturbations is the least.
More generally, this remains an open question, and testing this on different systems would be an interesting direction for future research.

There are significant areas for future work that remain outstanding. In principle, our method can be directly applied to a more realistic ephemeris model, with elliptical orbits for the primaries, external perturbations, and a lower bound on the impulse magnitude. Since our method is data-driven, any additional perturbations would be included in the mapping discovery. Applying our method to a more realistic ephemeris model would be an interesting next step to verify our strategy as a real-world station-keeping method.  
Another challenge is that the size of the augmentation perturbation needed to discover an accurate mapping is not known \textit{a priori}. Since ensemble learning provided insight into the compromise that must be made between the variance of the mapping coefficients and the variance in the training data, we suggest that an active-learning framework \cite{Fasel2021Ensemble-SINDy:Control} could be used to alter the sampled range of UPOs and perturbation size. This would direct the sampling of the underlying UPO structure of a system to regions of phase-space where the model uncertainty can be reduced the most.  
In the case of higher-dimensional systems, the number of polynomial library terms increases exponentially and it becomes increasingly difficult to discover a sparse mapping \cite{Bramburger2020Data-DrivenOrbits}. Introducing dimensionality reduction methods may be necessary to apply our method to such systems though
techniques such as \cite{Champion2019Data-drivenEquations} which utilises autoencoders to find a coordinate transformation in which the dynamics are sparse or \cite{Vlachas2022MultiscaleDynamics} that incorporate autoencoders with recurrent neural networks to evolve the latent space dynamics. Ensuring that these methods preserve key properties of chaotic systems may be necessary, as chaotic systems are especially sensitive to model error and preserved quantities not being conserved. Specialised autoencoder architectures that incorporate echo state networks \cite{Ozalp2025StabilitySpaces}, neural ODEs \cite{Zeng2022Data-drivenLearning}, and symmetries \cite{Zeng2021SymmetryDynamics} have shown promise in this regard.

The benefits of our sample-efficient and verifiable method for discovering accurate local Poincar\'e maps and subsequent low-energy control of UPOs may be realised in other fields. One potential application is the control of UPOs in fluids where there is increasing knowledge of these structures \cite{Budanur2017RelativeFlow,Singer1991ActiveConvection, Page2020SearchingDecomposition, Page2024RecurrentStructures}, which we can utilise in our sample-efficient method to avoid additional expensive data collection. A first step in experimental verification could be to test our method on a double pendulum system, a table-top analogue of the 3BP \cite{Kaheman2023SaddlePendulum}. Additionally, other applications include energy-efficient control of robotics systems \cite{Morimoto2005Poincare-Map-BasedWalking} and plasma physics, an area which SINDy has already shown promise in \cite{Bayon-Bujan2024Data-drivenPropulsion, Alves2022Data-drivenSimulations} and can be modelled using discrete mappings \cite{Burby2020FastFields}. 
The work in \cite{Morimoto2005Poincare-Map-BasedWalking} used a model-based reinforcement learning framework in which a representation of a Poincar\'e map was learned in order to stabilise the natural gait  of a bipedal robot. In this vein, unifying the recently developed SINDy-RL (reinforcement learning) framework, which demonstrated a significantly increased sample efficiency \cite{Zolman2024SINDy-RL:Learning}, with our method may further reduce the number of samples needed to learn an effective Poincar\'e map-based controller as part of the aforementioned active-learning framework. This is the subject of ongoing research and we hope to report on developments in upcoming papers.

\section*{Acknowledgements}
This work was supported by the EPSRC [grant number EP/W524323/1].

\bibliographystyle{ieeetr}
\bibliography{Refs}

\end{document}